\documentclass{article}
\usepackage{amsfonts}
\usepackage{amsmath,latexsym}

\author{A.L.Lukashov\thanks{Supported by the Austrian Science Fund FWF,
project-number P16390-N04.}
{}\thanks{Research of the first author was supported  partly  by grant of President of Russian Federation, grant
NSh-1295.2003.1} and F.Peherstorfer}
\title{Zeros of polynomials orthogonal on two arcs of the unit circle}

\newcommand{\inte}{\mathop{\rm int}\nolimits}
\newcommand{\abs}[1]{\lvert#1\rvert}
\begin{document}

\newpage
\pagenumbering{arabic}
\maketitle
  \begin{abstract}
In this paper we study polynomials $(P_n)$ which are hermitian 
orthogonal on two
arcs of the unit circle with respect to weight functions which have 
square root
singularities at the end points of the arcs, an arbitrary nonvanishing trigonometric
polynomial $\mathcal{A}$ in the denominator and possible point 
measures at the zeros of
$\mathcal{A}$. First we give an explicit representation of the 
orthogonal
polynomials $P_n$ in terms of elliptic functions. With the help of this
representation for sufficiently large $n$ the number of zeros of $P_n$ 
which are
in an $\varepsilon$-neighbourhood of each of the arcs are
determined. Finally it is shown that the accumulation points of the 
zeros of
$(P_n)$ which are not attracted to the support lie on a Jordan arc running within 
the unit
disk from one of the arcs to the other one. The accumulation points lie dense on 
the
Jordan arc if the harmonic measures of the arcs are irrational. If the 
harmonic
measures are rational then there is only a finite set of accumulation 
points on
the Jordan arc.

   \end{abstract}

\section{Introduction}
Let $d\leq\varphi_{1}<\varphi_{2}<\varphi_3<\varphi_4<d+2\pi$
and put $$E=[\varphi_1,\varphi_2]\cup[\varphi_3,\varphi_4]=E_1\cup E_2,$$
and let
$$\{z=e^{i\varphi}:\varphi\in E\}=\Gamma_E=\Gamma_{E_1}\cup\Gamma_{E_2}.$$
For 
\begin{math}
	n\in\mathbb{N}_0=\{0,1,2,\ldots\} 
\end{math} let  
$$\Pi_{n/2}=\left\{\sum_{k=0}^{[n/2]  }a_k\cos(\frac{n-2k}{2 }\varphi)+b_k
\sin(\frac{n-2k}{2 }\varphi):a_k,b_k\in\mathbb{R}   \right\} $$ 
denote the space of real trigonometric polynomials of (integer or half-integer) 
degree no more than $n/2.$  We say \( \mathcal{D}\in\Pi_{n/2}  \) is of exact degree 
$\partial\mathcal{D}=n/2,$ if $|a_0|+|b_0|\neq 0.$ By ${\cal 
R}\in\Pi_2$ we denote the trigonometric polynomial which vanishes at the 
endpoints of the two arcs, i.e.,
\begin{equation}
\label{eq-1.3}
 {\cal R}(\varphi)=\prod
_{k=1}^{4}\sin\frac{\varphi-\varphi_{k}}{2}
\end{equation}
and let
\begin{equation}
\label{eq-1.4} 	
{\cal R}(\varphi)={\cal V}(\varphi){\cal W}(\varphi)	
\end{equation} be an arbitrary splitting of  $\cal R$ with ${\cal V},{\cal W}\in\Pi_2.$ 

Loosely speaking we study polynomials which are orthogonal on the two arcs 
$\Gamma_E$ of the unit circle with respect to a distribution of the form 
$$\sqrt{\abs{\mathcal{W(\varphi)}}}/\mathcal{A(\varphi)}\sqrt{\abs{\mathcal{V(\varphi)}}}d\varphi+\mbox{ possible point measures at the zeros 
of }{\cal A}(\varphi),$$ where  ${\cal A}(\varphi)$ is a real trigonometric polynomial 
which has no zeros in $E$ and satisfies some other mild conditions, see 
\eqref{eq-1.5} below, also concerning the precise form of the point measures. (In 
fact even more general distributions including sign changing ones are 
considered). 

First we give an explicit representation of the orthogonal
polynomials in terms of elliptic functions and show how this 
representation can
be used also to obtain trigonometric polynomials minimal on two 
intervals with
respect to a weight function of the type 
$1/\sqrt{\abs{\mathcal{A}}}$. Then we
emphasize on the zeros of the orthogonal polynomials. Let
us recall that it is known by Fejer's Theorem on zeros of minimal 
polynomials \cite{fej}
that all zeros of $P_n$ lie in  the convex hull of $\Gamma _E$ (in fact,
strictly inside, by \cite{saf}) and that 
they are
attracted to the support up to a finite number (Widom's theorem \cite{wid1}). 
Furthermore it is known (see e.g. \cite[Th.5.2]{saf},\cite{statot}) that the 
zero distribution of $(P_n)$ converges weakly to the equilibrium distribution 
of $\Gamma_E,$ i.e.
\begin{equation}
\label{dirak}
\frac{1}{n}\sum_{j=1}^{n}\delta_{z_{j,n}}\xrightarrow[n\to\infty]{}\nu_{\Gamma_E},
\end{equation}
where $\delta_{z_{j,n}}$ denotes, as usual, the Dirac-delta measure at the 
point $z_{j,n}$ and $\nu_{\Gamma_E}$ the equilibrium measure of $\Gamma_E.$

Naturally we would like to know the precise number of the zeros of $P_n$ 
attracted to each of the two arcs
and what about the zeros which 
are not
attracted to the support. Concerning the first question we present a 
formula for
the precise number of zeros which are in an $\varepsilon$-neighbourhood 
of each of the arcs for sufficiently large $n.$ Then the behaviour of the accumulation points of the zeros of 
$(P_n)$
not attracted to the support is investigated. It is shown that they lie 
on an
open analytic arc with endpoints which are inner points of $\Gamma_{E_j},$ 
respectively, and can be given explicitely, see \eqref{curve} below.
Furthermore, the set of accumulation points is dense 
on this curve
if the harmonic measures of the arcs are irrational. If the harmonic 
measures are
rational then the set of accumulation points of zeros on the analytic arc is 
finite.
The last case  is known already \cite[Th.3.3]{pehste3} when one takes into 
consideration the known
fact that the reflection coefficients are pseudoperiodic if and only if the 
harmonic
measures are rational (see \cite{pehste4} and concerning pseudoperiodicity 
\cite[Th.1(a)]{luk1}).

Let us note that the behaviour of zeros of polynomials
orthogonal on the whole unit circle  is very different 
from that
one in the two arcs case. Indeed, it is well known in the case of the whole 
unit
circumference zeros need not be attracted to the support as the simple 
example
$P_n(z)=z^n$ shows. Let us mention also that  in the case of one arc, under the assumption
that the weight function is sufficiently nice, the reflection coefficients 
converge and thus there is always at 
most one point
(which can be deduced from \cite{pehste3}) to which zeros may be attracted if they are not attracted to 
the
support.

Using the fact that weight functions of the form
$\sqrt{\abs{\mathcal{W}}}/f\sqrt{\abs{\mathcal{V}}}$ on $E$ and 
zero
otherwise can be approximated well by weights 
$\sqrt{\abs{\mathcal{W}}}/\mathcal{A}_n\sqrt{\abs{\mathcal{V}}}$ treated in this
paper it can be shown using Tomcuk's asymptotic approach \cite{tom} (compare also 
\cite{wid2}) that the 
polynomials
orthogonal with respect to the above weights are asymptotically equal 
and that the
behaviour of the zeros is the same also, that is, is such as described 
in this
paper. This will be demonstrated in a forthcoming paper \cite{lukpeh}. At this 
point let us mention that asymptotic representations of polynomials orthogonal 
on two arcs of the unit circle can be obtained also from the very general and 
nice results of Widom \cite{wid2}. To extract the behaviour of the zeros of the 
orthogonal polynomials his results seem to be not explicit enough (compare 
\cite{apt} with this respect also).

Let us mention that by \cite[Th.2.1.3]{statot} there exists measures such that the set of 
accumulation
points of the zeros of $(P_n)$ is dense in the convex hull. For measures whose 
support is
the unit-circumference which have the property that the accumulation 
points of
the zeros are dense in $\abs{z}<1$, so-called Turan measures, see the 
recent
discussions in \cite{khr, sim}.
The results of this paper should be compared with the results on the 
zeros of polynomials
orthogonal on two intervals $[-1,a] \cup [b,1]$, $-1<a<b<1$, where a 
similar
behaviour of the zeros has been observed by the second author \cite{pehel2}
concerning the number of zeros in the intervals $[-1,a]$, $[b,1]$ and 
the
denseness of zeros in the gap $[a,b]$. In the meantime these results 
have been
extended to several intervals \cite{pehgen}, see also \cite[p.92]{sue1} for denseness 
results
under certain assumptions.

There is also a vaste literature dedicated to similar questions about zeros of nonhermitian orthogonal 
polynomials or more generally of denominators of Pad\'e-approximants. With this respect we refer  to the 
survey \cite{sue1} and the recent papers \cite{bus,lub}.

\section{Notations - Examples.}

Henceforth let
 ${\cal A}(\varphi)\in\Pi=\cup_{l=0}^{\infty}\Pi_l$ be  an
arbitrary  real trigonometric polynomial which has no zeros in $E,$ i.e.,
\begin{equation}
	\mathcal{A}(\varphi)\neq 0\mbox{ for }\varphi\in E,  
	\label{eq-1.5}
\end{equation} 
and thus $\mathcal{A} $  can be represented in the form
\begin{equation}
\mathcal{A}(\varphi)=c_{\mathcal{A}}\prod_{j=1}^{m^* }\left(\sin
\frac{\varphi-\xi_j}{2 }
 \right)^{ m_j},
	\label{eq-1.6}
\end{equation} 
where $m^*,m_j\in\mathbb{N} $  and where the $\xi_j$'s are distinct, lie
in $\mathbb{C}\setminus E $   and for $\xi_j\notin \mathbb{R} $  there
exists a 
$\xi_k=\overline{\xi_j} $  with $m_k=m_j.$

As announced in this paper we study polynomials $P_n$ orthogonal with respect 
to the functional
${\cal L}(\cdot;{\cal A,W},\lambda),$ i.e.
\begin{equation}
\label{1.20}
{\cal L}(z^{-k}P_n)=0,\quad\mbox{ for }  k=0,\ldots,n-1,
\end{equation}
where the functional is given as follows:
\begin{equation}
{\cal L}(h;{\cal A,W},\lambda):=\frac{1}{2\pi}\int\limits_E  h(e^{i\varphi}) 
f(\varphi;{\cal A,W})d\phi+\mathcal{G}(h;{\cal A,W},\lambda),	
	\label{eq-1.17}
\end{equation} 
with
\begin{equation}
f(\varphi,{\cal A},{\cal W})=
\genfrac{\{ }{.}{0pt}{}{\frac{{\cal W}(\varphi)}{{\cal A}(\varphi)r(\varphi)},
\qquad \varphi\in E,}{0, \qquad \varphi\notin E,}	
	\label{eq-1.7}
\end{equation} 
and
\begin{equation}
\mathcal{G} (h;{\cal A,W},\lambda)=\frac{1}{2}\sum_{j=1}^{m^*}
(1-\lambda_j)\sum_{\nu=0}^{m_j-1 } \mu_{j,\nu}(-1)^{\nu} \delta_{z_j}^{(\nu)}  \left(\frac{h(z)}{z } \right), 
		\label{eq-1.18}
\end{equation} 
where
\begin{equation}
		\frac{1}{r(\varphi) } :=\frac{(-1)^{j} }{\sqrt{|\mathcal{R}(\varphi)| }  } ,\quad 
	j=1,2;
	\label{eq-1.8}
\end{equation} 
the $\mu_{j,\nu} $ 's are certain complex numbers (for their exact description see 
\eqref{eq-2.6} below)depending on 
$\mathcal{A},\mathcal{W}  $  and $\mathcal{R} $, $z_j:=e^{i\xi_j} \in\mathbb{C} \setminus\Gamma_E, 
\delta_{z_j}^{\nu}(g):=(-1)^{\nu}g^{(\nu)}(z_j)/\nu!    $  and
\begin{equation}
	\lambda=(\lambda_1,\ldots,\lambda_{m^*}),\mbox{ with 
	}\lambda_j\in\{-1,1\}\mbox{ and such that 
	}\lambda_{j_1}=\lambda_{j_2}\mbox{  for } \xi_{j_1}=\bar{\xi}_{j_2}.         
	\label{eq-1.19}
\end{equation} 
The functional ${\cal L}(\cdot;{\cal A,W},\lambda)$ was introduced in 
\cite{pehste1} for an arbitrary number of arcs even. Naturally the functional need not be positive definite.
As usual we call a functional $\mathcal{L} $  positive-definite, if 
$\det(c_{j-k} )_{j,k=0}^n>0 $  for all $n\in\mathbb{N}_0, $  
where the moments $c_j$ are given by $c_j=\mathcal{L}(z^{-j}), j\in\mathbb{Z}.   $ 
Note that $\mathcal{L} $ is positive definite if $f$ has no sign change on $E$ 
and $m_j=1$ for all $j\in\{1,\ldots, m^*\}$ for which $\lambda_j=-1.$
For $\partial\mathcal{W} = \partial\mathcal{V}-2,\lambda\equiv 1$ and  
$\mathcal{L} $ positive definite we obtain weights  studied by Yu.Ya.Tomchuk 
\cite{tom}. For other studies of polynomials orthogonal with respect to  
$\mathcal{L}, $ see also \cite{gerjoh}.
If  $\mathcal{L} $ is not positive definite we may have higher 
orthogonality of $P_n,$ indeed we may have
$${\cal L}(z^{-k}P_n)=0,\quad\mbox{ for }  k=0,\ldots,n+\mu-1,\mu\in\mathbb{N}.$$
As we shall see (Corollary 2 to Theorem 1 below, compare also \cite[Th.1]{luk2}) this case of maximal 
orthogonality is of interest in describing rational trigonometric functions 
which  deviate least from zero on the two intervals $E_1\cup E_2.$ 

Let us give two examples, the first one related to the positive definite case 
and the second one related to higher orthogonality.

EXAMPLE 1. Suppose that $\mathcal{A} $  has only real zeros and exactly one 
simple zero in each interval 
$(\varphi_{2j},\varphi_{2j+1}),j=1,2;\varphi_5:=\varphi_1+2\pi,  $  
i.e.$$\mathcal{A}(\varphi)=\sin\left(\frac{\varphi-\xi_1}{2 } \right)
 \sin\left(\frac{\varphi-\xi_2}{2 } \right), $$ 
where $\xi_j\in(\varphi_{2j},\varphi_{2j+1}).  $  Then for the weight $f(\varphi;
\mathcal{A},1) $  the orthogonality condition (\ref{1.20})  takes the form, by 
inserting the explicit expressions for $\mu_{j,0} $ given in \eqref{eq-2.6}, 
\begin{equation}
\int_{E}^{ } \frac{e^{-ik\varphi}P_n(e^{i\varphi})  }{|\mathcal{A}(\varphi)|\sqrt{|\mathcal{R}(\varphi)| }   } d\varphi+\sum_{j=1}^{2 } 
\frac{(1-\lambda_j)}{2 } \frac{\sqrt{R(e^{i\xi_j}) } e^{-i(k+1)\xi_j}P_n(e^{i\xi_j})  }{i\left(\frac{d}{dz }A\right)(e^{i\xi_j})    } =0
	\label{eq-1.26}
\end{equation} 
for $k=0,\ldots,n-1,$  where 
$R(e^{i\varphi}):=e^{2i\varphi} \mathcal{R}(\varphi)  $  and 
$A(e^{i\varphi}):=e^{i\varphi} \mathcal{A}(\varphi)  .$  Recall that 
$\lambda_1,\lambda_2$ can be chosen arbitrary from $\{-1,1\} .$   Relation 
(\ref{eq-1.26}) represents an orthogonality relation for $P_n$  with respect 
to a positive measure $d\sigma$ which has mass points at those $e^{i\xi_j} $   
where $\lambda_j=-1.$ 

2.If there exists a $T$-polynomial $\mathfrak{T}_N $  on $E$  (see 
\cite{pehste2}), then it is orthogonal with respect to the sign-changing weight 
$f(\varphi;1,1),$   namely,
$$\mathcal{L}(z^{-k}\mathfrak{T}_N;1,1,1)=0\mbox{ for } k=0,\ldots,N,    $$ 
and, denoting by $\alpha,|\alpha|=1,$  the leading coefficient of 
$\mathfrak{T}_N, $  the trigonometric polynomial 
$\tau_N(\varphi)=e^{-i(N/2)\varphi}\mathfrak{T}_N(e^{i\varphi})   $  deviates 
least from zero on $E$ with respect to the sup-norm among all trigonometric 
polynomials of degree $N/2$ with leading coefficients $2\cos\psi$  and $2\sin\psi,$ 
where $\alpha=e^{-i\psi}. $

In the following we need the additional notations: let $\mathbb{P}_n$ denote 
the set of algebraic polynomials of degree $n,$ let $A(z)$ be the algebraic polynomial which is connected with ${\cal A}(\varphi)$  by the relation
\begin{equation}
\label{eq-el1}
A(e^{i\varphi})=e^{ia\varphi}{\cal A}(\varphi), 
\end{equation}
where $2a=2\partial\mathcal{A}=\sum_{j=1}^{m^*}m_j, $  i.e.
$$A(z)=c_A\prod_{j=1}^{m^*}(z-z_j)^{m_j},$$
with  $c_A\in\mathbb C,$ $z_j=e^{i\xi_j},j=1,\ldots,m^*,$
 and all $z_j$ are distinct, and for $|z_j|\neq 1$ there exists $k$ such that $z_k=1/\overline{ z_j}
 ,m_k=m_j.$ The polynomial $A$ coincides with its reciprocal polynomial
$$A^*(z)=z^{2a}\overline{A(1/\bar z)},$$ 
i.e. it is selfreciprocal. Furthermore, $R,V,W$ are algebraic polynomials of degrees 
$4,2v,2w$ correspondingly, which can be obtained from ${\cal R},{\cal V},{\cal W}$
in an analogous way to (\ref{eq-el1}),$w_j$ denotes the number of zeros of ${\cal 
W}$ on $[\varphi_{2j-1},\varphi_{2j}],j=1,2,$  
$$A_j(z)=\frac{A(z)}{(z-z_j)^{m_j}},$$
and we make the additional supposition $a-w+1\in\mathbb{ N}_0.$

Now we can describe more precisely the functional $\mathcal{G} (h;{\cal 
A,W},\lambda)$  from \eqref{eq-1.18}. Namely, 
\begin{equation}
	\mathcal{G} (h;{\cal A,W},\lambda)=\frac{1}{2 } \sum_{j=1}^{m^* }
\frac{1-\lambda_j}{(m_j-1)! } 
	\left(\frac{z^{a-w}Wh }{iA_j\sqrt{R}  } \right) ^{(m_j-1)} (z_j),
	\label{eq-2.6}
\end{equation} 
i.e.
$$\mu_{j,\nu}=\frac{1}{(m_j-1-\nu)! } 	\left(\frac{z^{a-w}Wh }{iA_j\sqrt{R}  } \right) ^{(m_j-1)} (z_j), $$ 
where here and everywhere later by $\sqrt{R} $   the branch on 
$\mathbb{C}\setminus\Gamma_E $ is denoted which satisfies
\begin{equation}
	\arg\sqrt{R(e^{i\varphi}) } =\arg(-e^{i\varphi} ),\quad 
	\varphi\in(\varphi_2,\varphi_3).
	\label{eq-el2}
\end{equation} 
In the case when $m_j=1,j=1,\ldots,m^*,$ the functional $\cal L$ is nothing else as 
the Stieltjes 
integral with respect to the measure with absolute continuous part
$f(\varphi,{\cal A},{\cal W})d\varphi$ and with possible addition of masses at the points $z_j$.

So the main objects of investigation are the polynomials $P_n$, which are orthogonal with respect
to the functional $\cal L$ in the sense of \eqref{1.20} .
We shall use the notation
$${\cal L}(z^{-k}P_n)=0,\quad k\in(0,n-1)$$
for (\ref{1.20}). But if it is known that  ${\cal L}(z^{-n}P_n)\neq
0$, then we shall write
${\cal L}(z^{-k}P_n)=0,k\in(0,n-1].$

The following conformal mapping of a certain rectangle in the complex plane to the exterior of
$\Gamma_E$ will play a crucial role in the statement of our results. Let
\begin{equation}
\label{eq-el4}
k^2=( e^{i\varphi_1}, e^{i\varphi_2}, e^{i\varphi_3}, e^{i\varphi_4})
\end{equation}
be the modulus of the exterior of $\Gamma_E$, where
\begin{equation}
\label{eq-el5}
(z_1,z_2,z_3,z_4):=\frac{z_4-z_1}{z_4-z_2}:\frac{z_3-z_1}{z_3-z_2}
\end{equation}
denotes the double relation between points $z_1,z_2,z_3,z_4.$
The modulus $k$ will be simultaneously the modulus of the Jacobian elliptic functions
  \[\operatorname{sn} z = \operatorname{sn}(z;k), \quad
                             \operatorname{cn} z = \operatorname{cn}(z;k) = \sqrt{1 
                             -\operatorname{sn}^2 z},\]
                             and
                             \[\operatorname{dn} z = \operatorname{dn}(z;k) = \sqrt{1 - k^2
                             \operatorname{sn}^2 z},\]
                             and let $K = K(k)$ be the complete elliptic integral of the first kind 
                             of
                             modulus $k$ defined by
                             \begin{equation}
                             \label{prozentneulabel}
                             K = K(k) = \int_0^1 \frac{dx}{\sqrt{(1-x^2)(1-k^2 x^2)}}.\end{equation}
                             As usual let
                             \[k' = \sqrt{1- k^2} \quad \text{and} \quad K' = K'(k')\]
                             denote the complementary modulus and the complete elliptic integral of 
                             the
                             first kind with respect to $k'$, respectively.
   Furthermore, in the following we assume  without loss of generality  that
   $$\varphi_1=2\pi-\varphi_4,$$
   since it can be satisfied after a
suitable turn 
   of the unit circle, and such a turn corresponds to a substitution of the 
   kind $z\rightarrow e^{i\psi}z. $ 
   
   Next let us construct the conformal mapping from the (partly open) rectangle
                             \[\Box = \{ u \in \mathbb{C} : -K < \Re u <
0,\; -K' < \Im u \le K'\}\] to the exterior of  $\Gamma_E.$ Since the conformal 
mapping $w(u)$ from \( \Box \) to the exterior of two disjoint intervals say $[-1,\alpha]\cup[\beta,1] 
,-1<\alpha<\beta<1$ is known to be (see \cite[p.~139]{akh}, 
\cite{fis})
\begin{equation}
	w(u)=\frac{\operatorname{sn}^2 u\operatorname{cn}^2 a+
	\operatorname{cn}^2 u\operatorname{sn}^2 a}{\operatorname{sn}^2 u-
	\operatorname{sn}^2 a}=\alpha+\frac{1-\alpha^2}{2(\operatorname{sn}^2 u-
	\operatorname{sn}^2 a)} 
	\label{omega18}
\end{equation}
where 
 \begin{equation}
                                     \alpha = 1 - 2 \operatorname{sn}^2 a,
                                     \label{alpha}
                             \end{equation}
							 and
 \[\beta  =  2 \operatorname{sn}^2 
                             (K+a)-1,\]							 
we obtain the desired mapping $z=\phi(u)$ easily by composition of $w$ with 
the  M\"obius map 
\begin{equation}
	z=\frac{w-i\tan\frac{\varphi_1}{2 } }{w+i\tan\frac{\varphi_1}{2 } }, 
	\label{moebius}
\end{equation} 
which maps the upper half plane to the interior of the unit disk and the 
intervals $[-1,\alpha]\cup[\beta,1]$ to the arcs  $\Gamma_{E_1}$ and $ 
\Gamma_{E_2}.$ Thus the function 
\begin{equation}
z=\phi(u)=\frac{2\operatorname{sn}^2u\sin\frac{\varphi_1}{2 }e^{i\varphi_2/2}+(\alpha-1)e^{i\varphi_1/2}   }{2\operatorname{sn}^2u\sin\frac{\varphi_1}{2 }
e^{-i\varphi_2/2}+(\alpha-1)e^{-i\varphi_1/2} } 	
	\label{map}
\end{equation} 
  where
 $$\alpha=-\tan\frac{\varphi_1}{2 }\cot\frac{\varphi_2}{2 } =
 1 - 2 \operatorname{sn}^2 a,  $$
$$\beta=-\tan\frac{\varphi_1}{2 }\cot\frac{\varphi_3}{2 }=2 \frac{\operatorname{cn}^2 a}{\operatorname{dn}^2 a}-1 $$
realizes that map.   It is an even elliptic function of order 2 with primitive periods $2 K$ 
                             and $2 i
                             K'$ and simple poles at $\pm \zeta$ in the period parallelogramm
                             \[{\mathcal P} = {\mathcal P}(k) = \{u \in
\mathbb{C} : -K \le \Re u < K,\; 
                             -K' < \Im u \le K'\},\]
where $\zeta\in\Box$  is defined by the relation
\[ 
\operatorname{sn}^2\zeta=\frac{\sin\frac{\varphi_1+\varphi_2}{2 }e^{i\frac{\varphi_2-\varphi_1}{2 } }  }{\sin\varphi_1\sin\frac{\varphi_2}{2 }  }.  \] 
The points
 \[z: e^{i\varphi_1}  \to e^{i\varphi_2} \to e^{i\varphi_3} 
 \to e^{-i\varphi_1}\to  e^{i\varphi_1}  \]
                             correspond under the map $\phi(u)$  to the points
                             \[u:  0 \to i K' \to -K + i K' \to -K \to 0 ,\]
                             and the upper and lower halves of the open rectangle, that is, $(-K,0) 
                             \times
                             (0,i K')$ and $(-K,0) \times (0,-i K')$, are mapped onto
the interior and exterior of the unit circumference, respectively. Furthermore, we 
                             need the
                             theta functions $H$ and $\theta$ defined by (see, for example, 
                             \cite{ref28})
                             \[H(z) = \delta_1 \left( \frac{z}{2 K}\right) = 2 \sum_{j=0}^\infty 
                             (-1)^j
                             q^{(j+\frac{1}{2})^2} \sin \frac{(2 j + 1) \pi}{2 K}z,\]
                             and
                             \[\theta (z) = \delta_4 \left( \frac{z}{2 K}\right) = 1 + 2 
                             \sum_{j=1}^\infty
                             (-1)^j q^{j^2} \cos \frac{j \pi}{K}z,\]
                             and related to each other by
                             \[H(z+i K') = i e^{- i \pi z /2 K}q^{-1/4} \theta (z),\]
                             where $q = e^{- \pi K'/K}$. Note that $H$ and $\theta$ 
                             is an odd and 
                             an even
                             function, respectively. Both are analytic at every point of the complex 
                             plane
                             and are quasi doubly-periodic functions, that is, 
                             they satisfy the relations
\begin{equation}
	H(z+2K)=-H(z),\qquad H(z+2iK')=-e^{-i\pi z/K}q^{-1}H(z),  
	\label{foreta}
\end{equation} 
\begin{equation}
	\theta(z+2K)=\theta(z),\qquad\theta(z+2iK')=-e^{-i\pi z/K}q^{-1}\theta(z).  
	\label{fortheta}
\end{equation} 
\newcounter{tmp}													 
\section{The basic results}
\newenvironment{alphenumerate}
{\begin{list}{\alph{tmp})}{\usecounter{tmp}}}
{\end{list}}

\newtheorem{lemma}{Lemma}
\begin{lemma}\label{lem.1}
	\begin{alphenumerate}
	\item
The  Green's function $g$ of $\bar{\mathbb C}\backslash \Gamma_E$
with respect to the point $c_0\in\bar{\mathbb C}\backslash \Gamma_E$, 
is given in terms of  Jacobian elliptic functions by the relation
\begin{equation}
\label{eq-el7}
g_{\bar{\mathbb C}\backslash \Gamma_E}(z,c_0)=\log \left|\frac{H(u+\bar{\gamma})}{H (u-\gamma)}\right|,
\end{equation}
where $z=\phi(u),\phi$   is given by \eqref{map} and $c_0=\phi(\gamma).$ In 
particular, $\infty=\phi(\zeta).$
\item The harmonic measure of  $\Gamma_{E_2}$ at $z=\infty$ is given by 
\begin{equation}
\label{eq-el6}
\omega_2(\infty)=-\frac{\Re\zeta}{K}.
\end{equation}
Note that $\omega_1(\infty)+\omega_2(\infty)=1.$
\item The capacity of  $\Gamma_E$ is as follows
\begin{equation}
\tau=\operatorname{cap}	(\Gamma_E)=\left|\frac{H(2i\Im\zeta)}{H(2\zeta)}\right|.
	\label{capac}
\end{equation} 
\end{alphenumerate}
\end{lemma}

Proof ad a).

The function $g(z,c_0)$ defined by  (\ref{eq-el7}) 
is harmonic on $\bar{\mathbb C}\backslash(\Gamma_E\cup\{c_0\})$ since it is a 
single-valued real part of the multi-valued analytic function
$$\log \frac{H(u+\bar{\gamma})}{H(u-\gamma)},$$
which follows by the facts that 
$$\frac{H(u+2iK'+\bar{\gamma})}{H(u+2iK'-\gamma)}: \frac{H(u+\bar{\gamma})}{H(u-\gamma)}=e^{-2i\pi\Re\gamma/K}$$
and that
$$\frac{H(u+2K+\bar{\gamma})}{H(u+2K-\gamma)}: \frac{H(u+\bar{\gamma})}{H(u-\gamma)}=1.$$
Furthermore, for $\gamma\neq\zeta\mbox{ (i.e. } c_0\neq\infty)$
\begin{equation}
\label{eq-el8}
\log\left|\frac{H(u+\bar{\gamma})}{H(u-\gamma)}\right|+\log|u-\gamma|=\log\left|\frac{H(u+\bar{\gamma})}{H(u-\gamma)}\cdot(u-\gamma)\right|
\end{equation}
is a bounded function in a neighbourhood of $\gamma$, hence 
$g(z,c_0)=-\log|z-c_0|+O(1),$ as $z\to c_0.$
Analogously $g(z,\infty)=\log|z|+O(1),$ as $z\to\infty.$
Moreover, for $\Re u=0,-K'\leq\Im u\leq K',$
\begin{equation}
\label{eq-el9}
\left|\frac{H(u+\bar{\gamma})}{H(u-\gamma)}\right|^2=\frac{H(u+\bar{\gamma})}{H(u-\gamma)}\frac{H(\bar u+\gamma)}{H(\bar u-\bar{\gamma})}=1,
\end{equation}
and analogously (\ref{eq-el9}) holds for $\Re u=-K,-K'\leq\Im u\leq K'.$ Hence 
(\ref{eq-el7}) is the Green's function.

Proof ad c).

From the definition of capacity it follows that
\begin{equation}
	\tau=e^{-\gamma},\gamma=\lim_{z\to\infty}g(z,\infty)-\log|z|.
	\label{proofcap}
	\end{equation}
Since by \eqref{map} and the Representation theorem for elliptic 
functions
$$|\phi(u)|=\left|\frac{H(u-\bar{\zeta})H(u+\bar{\zeta})}{H(u-\zeta)H(u+\zeta)}\right|$$
the relations \eqref{eq-el7} and \eqref{proofcap}  give the desired result.

The proof ad b) will be given in Theorem 2.     $\Box$

For another representation of the Green's function see \cite{pehste4}.

\par
The starting point of our investigations is the following characterization
(due to the second author and R.Steinbauer \cite{pehste1}) of the polynomials orthogonal
with respect to ${\cal L}(\cdot;{\cal A,W},\lambda)$
by a quadratic equation.
\par

\begin{lemma}
Let ${\cal L}(\cdot;{\cal A,W},\lambda)$ be given as above, let $a-w+1\in\mathbb 
N_0,$ and let $\mu\in\{0,1\}.$ Then for $n\geq a+1+v$ the following assertions are 
equivalent:
\begin{enumerate}
\item ${\cal L}(z^{-j}P_n;{\cal A,W},\lambda)=0$ for $j\in(0,n+\mu-1]$ 
\item  there exists a
 polynomial $Q_{n+2-2v}\in\mathbb P_{n+2-2v}$ and there exists a polynomial $g_{(n)}
 \in\mathbb P_{1-\mu}$ with
 $g_{(n)}(0)\neq 0$  such that
\begin{equation}
\label{eq-el10}
W(z)P_n^2(z)-V(z)Q_{n+2-2v}^2(z)=z^{n+p-(a+1-w)+\mu}A(z)g_{(n)}(z),
\end{equation}
where $p,0\leq p\leq 1,$ is the multiplicity of the zero of $P_n$ at $z=0$ and that
\begin{eqnarray}
\label{eq-el11}
\left(VQ_{n+2-2v}\right)^{(k)}  (z_j)=\lambda_j\left( \sqrt{R}P_n\right)^{(k)}(z_j),
 ,\\
\mbox{ for } k=0,\ldots,m_j-1;j=1,\ldots,2m^*,\nonumber\\
\left.\frac{VQ_{n+2-2v}}{\sqrt{R}P_n}\right|_{z=0}=1,\qquad 
V(0)Q^*_{n+2-2v}(0)=\sqrt{R(0)}P_n^*(0).\label{eq-el12}
\end{eqnarray}
\end{enumerate}
\end{lemma}

The basic theorem for what follows is the next one.
\newtheorem{theorem}{Theorem}
\begin{theorem}
Let $n\geq 
a+1+v$ and $\mu\in\{0,1\}.$ If ${\cal L}(z^{-j}P_n;{\cal A,W},\lambda)=0$ for 
$j\in(0,n+\mu-1]$  then the polynomials $P_n$ and $Q_{n+2-2v}$ from  Lemma 2 satisfy the following relations:
\begin{equation}\frac{2W(z)P_n^2(z)}{z^{n+p-(a+1-w)}A(z)g_{(n)}(z)}-1=
\label{eq-el13}
\frac{1}{2}\left(\Psi_n(\phi(u))+\Psi_n(-\phi(u))\right),
\end{equation}
and
\begin{equation}
\label{eq-el14}
\frac{2P_n(z)\sqrt{R(z)}Q_{n+2-2v}(z)}{z^{n+p-(a+1-w)}A(z)g_{(n)}(z)}=
\frac{1}{2}\left(\Psi_n(\phi(u))-\Psi_n(-\phi(u))\right),
\end{equation}
where
\begin{equation}
\label{eq-el15}
\Psi_n(\phi(u))=ce^{-i\pi m^{(n)}u/K} \left[\frac{H(u+\zeta)}{H(u-\zeta)}\right]^{n+w-p-a+1-\mu-\partial g_{(n)}} 
\end{equation}
$$\cdot\left[\frac{H(u+\bar{\zeta})}{H(u-\bar{\zeta})}\right]^{n+w-p-a-1+\mu} 
\prod_{j=1}^{m^*}\left[\frac{H(u+v_j)}{H(u-v_j)}\right]^{\lambda_jm_j}\left[\frac{H\left(u+b^{(n)}\right)}
{H(u-b^{(n)})}\right]^{\delta^{(n)}},$$
$b^{(n)}\in \Box,$   $m^{(n)} \in \mathbb Z,$ and $\delta^{(n)}\in\{-1,0,1\}$ 
$(\delta^{(n)}=0\Longleftrightarrow\partial g_{(n)}=0)$ are given by
the system of equations
\addtocounter{equation}{1}
\newcounter{z1}
\setcounter{z1}{\value{equation}}
$$\left\{ \begin{aligned}
m^{(n)} K'+(2-2\mu-\partial g_{(n)})\Im\zeta+\sum_{j=1}^{m^*}\lambda_jm_j\Im v_j+\delta^{(n)}\Im
b^{(n)}&=&0 &\qquad{\rm (\arabic{equation}) }\\
\addtocounter{equation}{1}
(2n+2w-2a-\partial g_{(n)}-2p)\Re\zeta+\sum_{j=1}^{m^*}\lambda_jm_j\Re 
v_j+\delta^{(n)}\Re b^{(n)}&=&-l_nK,&\qquad {\rm (\arabic{equation}) }
\end{aligned}\right.$$
where $l_n\in\mathbb N,$ and 
\begin{equation}
c=(-1)^{2a}F_w(e^{i\varphi_1}).	
	\label{cth}
\end{equation} 

Here  $$F_w(z)=\frac{W^2(z)-V^2(z)}{W^2(z)+V^2(z)}.$$
\end{theorem}

{\bf Proof.} 

Let us consider the function
\begin{equation}
\label{eq-el18}
\Psi_n(z)= \frac{W(z)\left(P_n(z)+\sqrt{\frac{V(z)}{W(z)}}Q_{n+2-2v}(z)\right)^2}{z^{n+p-(a+1-w)+\mu}A(z)g_{(n)}(z)},
\end{equation}
where $Q_{n+2-2v}(z)$ and $g_{(n)}(z)$ are the polynomials from Lemma 2. The function $\Psi_n$ is meromorphic on the Riemann surface $\cal S$ of the function $\omega=\sqrt{R(z)}$ (since it is a rational function of the variables $\omega,z$). The Riemann surface $\cal S$ is a compact Riemann surface of genus 1, and the mapping $I(z,\omega)\longrightarrow(z,-\omega)$ changes sheets of $\cal S$.

The function 
\begin{equation}
\label{eq-el19}
\Psi_{1,n}(z)= \frac{W(z)\left(P_n(z)-\sqrt{\frac{V(z)}{W(z)}}Q_{n+2-2v}(z)\right)^2}{z^{n+p-(a+1-w)+\mu}A(z)g_{(n)}(z)}
\end{equation}
corresponds to the function $\Psi_n$ under the map $I.$

Applying the map $z=\phi(u)$, one obtains two functions $\Psi_n(\phi(u))$ and $\Psi_{1,n}(\phi(u))$
which are well-defined on the rectangle $\Box.$ One can extend them onto the period parallelogramm
$${\cal P}={\cal P}(k)=\{u\in\mathbb C:-K\leq\Re u<K,-K'<\Im u\leq K'\}$$
by 
$$\Psi_n(\phi(-u))\stackrel{def}{=}\Psi_{1,n}(\phi(u))\mbox{ and }
\Psi_{1,n}(\phi(-u))\stackrel{def}{=}\Psi_n(\phi(u)).$$ Then it is possible 
to extend them onto the whole plane $\mathbb C$ by the double-periodicity with respect to $2K$ and $2K'.$ Since both functions $\Psi_n$ and
$\Psi_{1,n}$ are rational functions of the variables $\omega,z$, they are meromorphic on the surface $\cal S$. 
It is well known that the Jacobian elliptic functions uniformize the surface $\cal S$,
hence the functions $\Psi_n(\phi(u))$ and $\Psi_{1,n}(\phi(u))$ are elliptic.

Let us determine all zeros and poles of $\Psi_n(\phi(u)).$ First, we conclude from 
(\ref{eq-el18}),(\ref{eq-el19}) and (\ref{eq-el10}) that
\begin{equation}
\label{eq-el20}
\Psi_n(\phi(u))\Psi_{1,n}(\phi(u))\equiv 1,
\end{equation}   
hence if $u$ is a zero of $\Psi_n (\phi(u))$ 
then $-u$ is a pole of $\Psi_{1,n} (\phi(u)),$ and vice versa. Now from 
(\ref{eq-el18}) and (\ref{eq-el12})
\begin{itemize}
\item[(i)] $u=\zeta $ (which corresponds to $x=\infty )$ is a pole of multiplicity $n+w-p-a+1-\mu-\partial g_{(n)}$
of $\Psi_n(\phi(u))$ 
\end{itemize} 
and by (\ref{eq-el20})
\begin{itemize}
\item[(ii)] $u=-\zeta $ is a zero of multiplicity
$n+w-p-a+1-\mu-\partial g_{(n)}$
of $\Psi_n(\phi(u)).$ 
\end{itemize}
Moreover, by (\ref{eq-el18}) and (\ref{eq-el12}) 
\begin{itemize}
\item[(iii)] $u=\bar{\zeta}$ is a pole of $\Psi_n(\phi(u))$ of multiplicity $n-p-(a+1-w)+\mu,$
\end{itemize}
and by (\ref{eq-el20})
\begin{itemize}
\item[(iv)] $u=-\bar{\zeta}$ is a zero of $\Psi_n(\phi(u))$ of multiplicity $n-p-(a+1-w)+\mu.$
\end{itemize}
From (\ref{eq-el18}),(\ref{eq-el20}) and (\ref{eq-el11}) it follows that
\begin{itemize}
\item[(v)] $u=v_j$ is a zero (pole) of  multiplicity $m_j$ of $\Psi_n(\phi (u)),$ if $\lambda_j=-1  (+1),j=1,\ldots,m^*;$
\item[(vi)] $u=-v_j$ is a zero (pole) of  multiplicity $m_j$ of $\Psi_n(\phi (u)),$ if $\lambda_j=+1 (-1),j=1,\ldots,m^*.$
\end{itemize}
Finally, for $\partial g_{(n)}=1,$
\begin{itemize}
\item[(vii)] $u=b^{(n)}$ is a zero (pole) of $\Psi_n(\phi(u)) $ if $\delta^{(n)}=-1(+1),$
\item[(viii)] $u=-b^{(n)}$ is a pole (zero) of $\Psi_n(\phi(u))$ if $\delta^{(n)}=-1(+1).$
\end{itemize}
Here $b^{(n)}\in\Box$ and $\delta^{(n)}\in\{-1,1\}$ are defined by
\begin{equation}
\label{eq-el21}
V(\phi(b^{(n)}))Q_{n+2-2v}(\phi(b^{(n)}))=\delta^{(n)}\sqrt{R(\phi(b^{(n)}))}P_n(\phi(b^{(n)})).
\end{equation}
Summing up (i)-(viii) we get by the  Representation theorem
for elliptic functions in terms of theta functions (see, for example,
 \cite[p.54]{akh}) that $\Psi_n(\phi(u))$ has a representation
of the form
\begin{equation}
\label{eq-el22}
\Psi_n(\phi(u))=c^{(n)}  
\left[\frac{H(u+\zeta)}{H(u-\zeta)}\right]^{n+w-p-a-\mu-\partial g_{(n)}} \left[\frac{H(u+\bar{\zeta})}{H(u-\bar{\zeta})}\right]^{n+w-p-a-1+\mu}  
\end{equation}
$$\cdot\frac{H(u+\zeta)}{H(u-\tilde{\zeta}^{(n)} ) } \prod_{j=1}^{m^*}\left[\frac{H(u+v_j)}{H(u-v_j)}\right]^{\lambda_jm_j}
\left[\frac{H\left(u+b^{(n)}\right)}{H(u- b^{(n)})}\right]^{\delta^{(n)}},$$
where $\tilde{\zeta} ^{(n)}=\zeta-2l^{(n)} K-2m^{(n)} iK'; l^{(n)} ,m^{(n)}\in\mathbb Z.$

With the help of
\eqref{foreta} 
one obtains from (\ref{eq-el22}) the required representation (\ref{eq-el15}) up 
to the multiplicative constant $c$.

The formulas (\ref{eq-el13}),(\ref{eq-el14}) together with
\begin{equation}
1+\frac{2V(z)Q_{n+2-2v}^2(z)}{z^{n+p-(a+1-w)}A(z)g_{(n)}(z)}=\frac{1}{2}\left(\Psi_n(\phi(u))-\Psi_n(-\phi(u))\right),	
	\label{eq-el13pr}
\end{equation}
needed in the following and which is just a rewriting of (\ref{eq-el14})
follow from 
(\ref{eq-el18}),(\ref{eq-el19}) and (\ref{eq-el10}).

Writing down the condition of ellipticity $\Psi_n(\phi(u+2iK'))=\Psi_n(\phi(u))$ 
for $\Psi_n$ from \eqref{eq-el22} 
gives (\arabic{z1}) and 
\addtocounter{z1}{1}
(\arabic{z1}).

To compute the constant $c$, put $u=0$ in (\ref{eq-el22}). Then
\begin{equation}
\label{eq-el23}
\Psi_n(\phi(0))=c^{(n)} (-1)^{l^{(n)}+m^{(n)}}q^{(m^{(n)})^2}  e^{\pi i m^{(n)}\zeta/K}(-1)^{2a}  .
\end{equation}
On the other hand, by (\ref{eq-el13})
\begin{equation}
\label{eq-el24}
\Psi_n(\phi(0))= \Psi_n(e^{i\varphi_1})=F_w(e^{i\varphi_1}),
\end{equation}
hence by (\ref{eq-el23}) and (\ref{eq-el24}) equality (\ref{cth}) follows
with
$$c=c^{(n)}(-1)^{l^{(n)}+m^{(n)}}q^{(m^{(n)})^2}  e^{\pi i
m^{(n)}\zeta/K}.  $$    

The case $\partial g_{(n)}=0$ is considered in an analogous way.     

Let us give another representation for $m^{(n)}$. For that 
reason let us put $u=-K$ in \eqref{eq-el15}  and \eqref{eq-el13}.  Then 
$$\Psi_n(\phi(-K))=c(-1)^{m^{(n)} } $$
and
$$\Psi_n(\phi(-K))=F_w(e^{i\varphi_4}).$$
So,
\begin{equation}
\label{eq-el25}
(-1)^{m^{(n)}}=(-1)^{2a}F_w(e^{i\varphi_1})F_w(e^{i\varphi_4}),
\end{equation}
and therefore $m^{(n)}$ is even (odd) for all $n\geq 
a+1+v$ simultaneously.      $\Box$

\newtheorem{corollary}{Corollary}

\begin{corollary}\label{cor-1}
 Let $n\geq 
a+1+v$.	If the functional $\mathcal{L} $ is positive definite then the monic polynomials $P_n$ orthogonal with respect to  
$\mathcal{L} $ have a representation of the form:
\begin{equation}
	P_n(\phi(u))=\frac{1}{2}\left(\Omega_n(u)+\Omega_n(-u)\right), 
	\label{eq-cor1}
\end{equation}  
where
\begin{equation}
	\Omega_n(u)=C_{\Omega,n}\left(\frac{H(u+\bar{\zeta})}{H(u-\zeta)} \right)^n
	\frac{H(u+\delta^{(n)}b^{(n)})}{H(u+\bar{\zeta})}
		\label{eq-cor2}
\end{equation} 
$$\cdot\frac{\left(H(u+\bar{\zeta})H(u+\zeta)\right)^{w-a}}{e^{i\pi k^{(n)}u/K} }
\frac{\prod_{j=1}^{m^*}H^{m_j\left(\frac{1+\lambda_j}{2  }\right) } (u+v_j)
H^{m_j\left(\frac{1-\lambda_j}{2  }\right) } (u-v_j)}{\prod_{j=1}^{2w}H(u-u_j) },$$
$k^{(n)}=(m^{(n)}-\#\{u_j:\Im u_j=K'\})/2$ and
\begin{equation}
C_{\Omega,n}=2e^{i(\varphi_1+4\phi)n}\frac{H^n(2i\Im\zeta)}{H^n(2\zeta) }
\frac{H(2\Re\zeta)}{H(\zeta+\delta^{(n)}b^{(n)})}
\frac{e^{i\pi k^{(n)}\zeta/K}}{\left(H(2\Re\zeta)H(2\zeta)\right)^{w-a}}	
	\label{eq-cor3}
\end{equation} 
$$\cdot\frac{\prod_{j=1}^{2w}H(\zeta-u_j)}{\prod_{j=1}^{m^*}
H^{m_j\left(\frac{1+\lambda_j}{2}\right)}(\zeta+v_j)
H^{m_j\left(\frac{1-\lambda_j}{2}\right)}(\zeta-v_j)},$$
$\phi=\arg H(\zeta).$

Furthermore
$b^{(n)}\in \Box,$   $m^{(n)} \in \mathbb Z,$ and $\delta^{(n)}\in\{-1,1\}$ 
 are given uniquely by
the system of equations 
\addtocounter{equation}{1}
\newcounter{z2}
\setcounter{z2}{\value{equation}}
$$\left\{ \begin{aligned}
m^{(n)} K'+\Im\zeta+\sum_{j\in J}\lambda_jm_j\Im v_j+\delta^{(n)}\Im
b^{(n)}&=&0\qquad &\qquad \qquad{\rm (\arabic{equation}) }\\
\addtocounter{equation}{1}
(2n+2w-2a-1)\Re\zeta+\sum_{j=1}^{m^*}\lambda_jm_j\Re 
v_j+\delta^{(n)}\Re b^{(n)}&=&-l_nK,&\qquad \qquad{\rm (\arabic{equation}) }
\end{aligned}\right.$$
where $l_n-w_2$ is even (recall that $w_2$ is the number of zeros of $\mathcal{W}$ 
on $E_2).$
Moreover the polynomials $Q_{n+2-2v}$ can be represented as
\begin{equation}
Q_{n+2-2v}(\phi(u))=\frac{1}{2}\left(\Omega_n(u)-\Omega_n(-u)\right)
\sqrt{\frac{W(\phi(u))}{V(\phi(u))}}.	
	\label{eq-cor4}
\end{equation}
\end{corollary}

Proof.

Let us define the function $\Omega_n(u),u\in\Box,$ as
\begin{equation}
\Omega_n(u)=P_n(\phi(u))+\sqrt{\frac{V(\phi(u))}{W(\phi(u))}}Q_{n+2-2v}(\phi(u)),	
	\label{eq-cor5}
\end{equation} 
where the polynomial $Q_{n+2-2v}$ is given by \eqref{eq-el10} . Since the substitution $u\to -u$ corresponds to the change of the 
branch of $\sqrt{\frac{V(\phi(u))}{W(\phi(u))}},$ we have
\begin{equation}
\Omega_n(-u)=P_n(\phi(u))-\sqrt{\frac{V(\phi(u))}{W(\phi(u))}}Q_{n+2-2v}(\phi(u)).	
	\label{eq-cor6}
\end{equation} 
Now formulas \eqref{eq-cor1},\eqref{eq-cor4} follow immediately from 
\eqref{eq-cor5},\eqref{eq-cor6}.

 Let us prove the representation \eqref{eq-cor2}. From 
\eqref{eq-el18} and \eqref{eq-cor5} it follows 
\begin{equation}
\Psi_n(\phi(u))=
\frac{W(\phi(u))\Omega_n^2(u)}{(\phi(u))^{n-a-1+w}A(\phi(u))g_{(n)}(\phi(u)) }.	
	\label{eq-cor7}
\end{equation} 
Further, applying the Representation theorem for elliptic functions one gets 
\begin{align}
	A(\phi(u))= &\mbox{ const } 
		\prod_{j=1}^{m^* } \left(\frac{H(u-v_j)H(u+v_j)}{H(u-\zeta)H(u+\zeta) } \right)^{m_j} ,  
	\label{eq-cor8}  \\
	  g_{(n)}(\phi(u))=  &\mbox{ const } 
		\frac{H(u-b^{(n)} )H(u+b^{(n)} )}{H(u-\zeta)H(u+\zeta) } ,  
	\label{eq-cor9}  \\
	W(\phi(u))= & \mbox{ const } 
		\prod_{j=1}^{2w} \frac{H(u-u_j)H(u+u_j)}{H(u-\zeta)H(u+\zeta) },  
		\label{eq-cor10}  \\
	\phi(u)= &\mbox{ const } 
		\frac{H(u-\bar{\zeta} )H(u+\bar{\zeta})}{H(u-\zeta)H(u+\zeta) }.  
	 	\label{eq-cor11}
\end{align} 
Substituting \eqref{eq-el18},\eqref{eq-cor8}-\eqref{eq-cor11} into 
\eqref{eq-cor7} gives
$$\Omega_n^2(u)=\mbox{ const 
} \left[\frac{H(u+\bar{\zeta}) }{H(u-\zeta) } \right]^{2n}
\frac{H^{1-\delta^{(n)} }(u-b^{(n)})  H^{1+\delta^{(n)} }(u+b^{(n)})  }{e^{i\pi 
m^{(n)}u/K  }  }$$
$$\cdot\frac{\left(H(u+\bar{\zeta})H(u+\zeta) \right)^{2w-2a}  }{H^2(u+\bar{\zeta}) }
\frac{\prod_{j=1}^{m^* }H^{m_j(1+\lambda_j)}(u+v_j)  H^{m_j(1-\lambda_j)}(u-v_j) }
{\prod_{j=1}^{2w}H(u-u_j)H(u+u_j) } =$$
\begin{equation}
	\mbox{ const 
} \left[\frac{H(u+\bar{\zeta}) }{H(u-\zeta) } \right]^{2n}
\frac{H^2(u+\delta^{(n)}b^{(n)}) }{e^{i\pi 
\tilde{m}^{(n)}u/K  }  }
\label{omegasqu}
\end{equation}
$$\cdot\frac{\left(H(u+\bar{\zeta})H(u+\zeta) \right)^{2w-2a}  }{H^2(u+\bar{\zeta}) }
\frac{\prod_{j=1}^{m^* }H^{m_j(1+\lambda_j)}(u+v_j)  H^{m_j(1-\lambda_j)}(u-v_j) }
{\prod_{j=1}^{2w}H^2(u-u_j)},$$
where $\tilde{m}^{(n)}=m^{(n)}-\#\{u_j:\Im u_j=K'\}.$ By the ellipticity of 
$\Omega_n$ and by \eqref{omegasqu} $\tilde{m}^{(n)}$ is even, 
$\tilde{m}^{(n)}=:2k^{(n)}$ which implies \eqref{eq-cor2} up to 
a constant multiplier.

To get \eqref{eq-cor3} one needs to take into account the equality 
$$1=\lim_{z\to \infty}\frac{P_n(z)}{z^n}=\lim_{u\to\zeta}
\frac{\frac{1}{2}\left(\Omega_n(u)+\Omega_n(-u)\right)}{\left(\mbox{ const } 
		\frac{H(u-\bar{\zeta} )H(u+\bar{\zeta})}{H(u-\zeta)H(u+\zeta) }\right)^n   }$$
and that the constant in \eqref{eq-cor11} can be easily determined from 
$\phi(0)=e^{i\varphi_1}.$

From the definition \eqref{eq-cor5} of the function $\Omega_n(u)$ it follows 
that it is a meromorphic function of $z=\phi(u)$ on the Riemann surface 
$\mathcal{S}$ of the function $\sqrt{R(z)}.$ Then, as it is known, $\Omega_n(u)$ 
is elliptic. Writing down the conditions of ellipticity for $\Omega_n(u)$ gives,
with the help of 
\eqref{foreta}, the relations (\arabic{z2}),
\addtocounter{z2}{1}
(\arabic{z2}) and that $l_n-w_2$ is even.

Conversely, let the conditions
\addtocounter{z2}{-1}
(\arabic{z2}),
\addtocounter{z2}{1}
(\arabic{z2}) with even $l_n-w_2$ be satisfied. Then the function $\Omega_n(u),$
defined by \eqref{eq-cor2}, is elliptic. Hence it can be represented as 
$\frac{p+\sqrt{R}q}{r},$ where $p,q,r$ are polynomials. From \eqref{eq-cor2} it 
follows that $\Omega_n$ as a function of $z$  has finite poles only at the zeros 
of $W$, and of the same order, hence $r(z)=W(z).$ Multiplying $\Omega_n(u)$ by 
$\Omega_n(-u)$ gives
$$\frac{p^2(z)-R(z)q^2(z)}{W^2(z)}=C_{\Omega,n}^2
\left(\frac{H(u-\bar{\zeta} )H(u+\bar{\zeta})}{H(u-\zeta)H(u+\zeta) }\right)^n   
\frac{H(u-\delta^{(n)}b^{(n)} )H(u+\delta^{(n)}b^{(n)} )}
{H(u-\bar{ \zeta})H(u+\bar{ \zeta}) }$$
$$\cdot\left(H(u-\bar{\zeta} )H(u+\bar{\zeta})H(u-\zeta)H(u+\zeta)\right) ^{w-a} 
\frac{	\prod_{j=1}^{m^* }H^{m_j} (u-v_j)H^{m_j}(u+v_j)}
{\prod_{j=1}^{2w}  H(u-u_j)H(u+u_j) },$$
what is equal by \eqref{eq-cor8}-\eqref{eq-cor11}  to
$$\mbox{ const 
}\frac{z^{n-a-1+w}  A(z)g_{(n)}(z) }{W(z) }=\frac{p^2(z)-R(z)q^2(z)}{W^2(z)},$$
i.e.,
\begin{equation}
	\frac{p^2(z)}{W(z) }-V(z)q^2(z)=\mbox{ const }z^{n-a-1+w}  A(z)g_{(n)}(z). 
	\label{eq-cor12}
\end{equation} 
Hence $p(z)=\tilde{P}(z)W(z),$ where $\tilde{P}(z)$ is a polynomial. Finally we get
$$\Omega_n(u)=\tilde{P}(\phi(u))+\sqrt{\frac{V(\phi(u))}{W(\phi(u))}}q(\phi(u)).$$
By \eqref{eq-cor2} $\Omega_n(u)$ has a pole of multiplicity $n$ at $u=\zeta,$ 
and $\Omega_n(-u)$ has a pole of multiplicity $w-a<n$ at $u=\zeta,$ hence 
$\tilde{P}(\phi(u))=\frac{1}{2}\left(\Omega_n(u)+\Omega_n(-u)\right)$ is a polynomial of 
degree $n$. Comparing the degrees in \eqref{eq-cor12}  gives  that 
$q$ is a polynomial of degree $n+2-2v.$ Putting $\tilde{P}_n :=\tilde{P}$ and $\tilde{Q}_{n+2-2v}:=q$  \eqref{eq-cor12}  becomes
$$  W(z)\tilde{P}_n^2(z)-V(z)\tilde{Q}_{n+2-2v}^2(z)=z^{n-(a+1-w)}A(z)\tilde{g}_{(n)}(z),$$
where $\tilde{g}_{(n)}$ is a polynomial of degree 1. Using representations  
\eqref{eq-cor5},\eqref{eq-cor6} and \eqref{eq-cor2} one gets after the
substitutions $u=\pm v_j, u=\pm\zeta$ the equalities 
\begin{eqnarray*}
\left(V\tilde{Q}_{n+2-2v}\right)^{(k)}  (z_j)=\lambda_j\left( \sqrt{R}\tilde{P}_n\right)^{(k)}(z_j),
 ,\\
\mbox{ for } k=0,\ldots,m_j-1;j=1,\ldots,2m^*,\\
\left.\frac{V\tilde{Q}_{n+2-2v}}{\sqrt{R}\tilde{P}_n}\right|_{z=0}=1,\qquad 
V(0)\tilde{Q}^*_{n+2-2v}(0)=\sqrt{R(0)}\tilde{P}_n^*(0).
\end{eqnarray*}
By Lemma 2  $\tilde{P}_n$ is 
orthogonal with respect to ${\cal L}.$ Since for a positive definite functional
${\cal L}$ the orthogonal polynomials are unique up to a multiplicative 
constant, the uniqueness of the solutions of the systems
\addtocounter{z2}{-1}
(\arabic{z2}) and
\addtocounter{z2}{1}
(\arabic{z2})  is proved.   $\Box$

\begin{corollary}
	Let $\mathcal{A}(\varphi)$ be a trigonometric polynomial 
	$\mathcal{A}\in\Pi_a,$ of the form \eqref{eq-1.6} which is positive on $E$ and 
	let $\nu\in\mathbb{N}/2$ be such that $\nu>a.$ Then the following assertions are 
	equivalent.
	\begin{alphenumerate}
		\item There exists a trigonometric polynomial 
		$$\tau_{\nu}(\varphi)=A\cos \nu\varphi+B\sin \nu\varphi+\ldots$$
		with $A,B,\in\mathbb{R},A^2+B^2\neq 0$ such that 
		\begin{equation}
			\max_{\varphi\in E}\left|\frac{\tau_{\nu}¥(\varphi)}
			{\sqrt{\mathcal{A}(\varphi)}}\right|=\min_{b_i,c_i\in\mathbb{R}}\left|
			\frac{A\cos \nu\varphi+B\sin \nu\varphi+b_1\cos(\nu-1)\varphi}
			{\sqrt{\mathcal{A}(\varphi)}}+\right.
			\label{eq-trig1}
		\end{equation}
		$$\left.\frac{c_1\sin (\nu-1)\varphi+\ldots+b_{[\nu]}
		\cos(\nu-2\left[\frac{\nu}{2}\right])\varphi+c_{[\nu]}
				\sin(\nu-2\left[\frac{\nu}{2}\right])\varphi}{\sqrt{\mathcal{A}(\varphi)}}\right|$$
	and all boundary points of $E$ are extremal points with
	\begin{equation}
		\frac{\tau_{\nu}(\varphi_{2j})}{\sqrt{\mathcal{A}(\varphi_{2j})}}=
		\frac{\tau_{\nu}(\varphi_{2j+1})}{\sqrt{\mathcal{A}(\varphi_{2j+1})}},\quad j=1,2.
		\label{eq-trig2}
	\end{equation}

	\item There exists a real trigonometric polynomial $\sigma_{\nu-1}\in\Pi_{\nu-1}$ 
	and a real constant $M_{\nu}$ such that
	\begin{equation}
		\tau_{\nu}^2(\varphi)-\mathcal{R}(\varphi)\sigma_{\nu-1}^2(\varphi)=
		M_{\nu}^2
		\mathcal{A}(\varphi)
		\label{eq-trig4}
	\end{equation}
	with
	\begin{equation}
		\tau_{\nu}(\xi_j)=(\sqrt{\mathcal{R}}\sigma_{\nu-1})(\xi_j),\quad j=1,\ldots,m^*.
		\label{eq-trig5}
	\end{equation}

        \item $\mathcal{L}(z^{-j}P_n;\mathcal{A},1,1)=0$ for $j\in (0,n],$ where

         $P_n(e^{i\varphi}):=e^{i\nu\varphi}\tau_{\nu}(\varphi),n=2\nu.$

	\item For some $l_{\nu}\in \mathbb{N}$ 
	\begin{equation}
		(4\nu-2a)\Re\zeta+\sum_{j=1}^{m^*}m_j\Re v_j=-l_{\nu}K
		\label{eq-trig3}
	\end{equation}
	holds.
	\end{alphenumerate}
	If any of those assertions holds then the minimal polynomial $\tau_{\nu}(\varphi)$ 
	is given by the formula
	\begin{equation}
	\tau_{\nu}(\varphi)=\frac{M_{\nu}}{2}\left(F_{2\nu}(u)+F_{2\nu}(-u)\right)
	e^{-i\nu\varphi},	
		\label{eq-trig6}
	\end{equation}
	$\exp(i\varphi)=\phi(u),$
	\begin{equation}
		F_{2\nu}(u)=\varepsilon_{\nu}e^{i\pi\mu/K}\left(\frac{H(u+\bar{\zeta})}{H(u-\zeta)}
		\right)^{2\nu}\frac{\prod_{j=1}^{m^*}H^{m_j}(u+v_j)}{\left(H(u+\zeta)
		H(u+\bar{\zeta})\right)^{a}},
		\label{eq-trig7}
	\end{equation}
	$|\varepsilon_{\nu}|=1,m=\frac{1}{2}\sum_{j=1}^{m^*}m_j\Im v_j\in\mathbb{N}.$
	\end{corollary}
	{\bf Proof.} The assertion is proved for the more general case of several arcs 
	(with expressions in terms of automorphic functions) in \cite[Th.1]{luk2}. We 
	give a proof also here for the sake of completeness.
	
	(a)$\Rightarrow$(b) Since
	$\{\frac{1}{\sqrt{\mathcal{A}(\varphi)}},
	\frac{\sin \varphi}{\sqrt{\mathcal{A}(\varphi)}},\ldots,
	\frac{\cos(\nu-1)\varphi}{\sqrt{\mathcal{A}(\varphi)}}\}$
	(for an integer $\nu$) and
\newline	$\{\frac{\sin\varphi/2}{\sqrt{\mathcal{A}(\varphi)}},
	\frac{\cos\varphi/2}{\sqrt{\mathcal{A}(\varphi)}}\ldots,
	\frac{\cos(\nu-1)\varphi}{\sqrt{\mathcal{A}(\varphi)}}\}$
	(for a half-integer $\nu$) are Chebyshev systems on $E$ by the Chebyshev 
	Alternation theorem we get that $\tau_{\nu}/\sqrt{\mathcal{A}}$ has at least 
	$2\nu-2$ alternation points  $\psi_j$ in the interior of $E.$ Put
	\begin{equation}
		\sigma_{\nu-1}(\varphi)=c\prod_{j=1}^{2\nu-2}\sin((\varphi-\psi_j)/2)
		\label{eq-corner1}
	\end{equation}
	and note that $\tau_{\nu}^2/\mathcal{A}-M_{\nu}^2$ has a double zero at any point 
	$\psi_j,j=1,\ldots,2\nu-2,$ and because of \eqref{eq-trig2} has a simple 
	zero at any zero of $\mathcal{R}(\varphi).$ Hence, for  
	a suitable constant $c$ in \eqref{eq-corner1},
	$$\frac{\tau_{\nu}^2(\varphi)}{\mathcal{A}(\varphi)}-M_{\nu}^2=
	\frac{\mathcal{R}(\varphi)
	\sigma_{\nu-1}^2(\varphi)}{\mathcal{A}(\varphi)},$$
	i.e. \eqref{eq-trig4} is proved.
	
	Furthermore, it follows from \eqref{eq-trig4} that in $[d,d+2\pi)\setminus E$ 
	the inequality
	$$\frac{|\tau_{\nu}(\varphi)|}{\sqrt{\mathcal{A}(\varphi)}}>M_{\nu}$$
	holds, hence the function
	\begin{equation}
      \label{trigf}
\frac{\tau_{\nu}(\varphi)+\sqrt{\tau_{\nu}^2(\varphi)-M_{\nu}^2\mathcal{A}
	(\varphi)}}
	{M_{\nu}\sqrt{A(\varphi)}}=:\frac{\mathcal{F}_{\nu}(\varphi)}{\sqrt{A(\varphi)}}
\end{equation}
	has on $E$ modulus 1 and on $[d,d+2\pi)\setminus E$ modulus greater  than 1. Furthermore, the function
	$$\frac{P_{2\nu}(z)+\sqrt{P_{2\nu}^2(z)-M_{\nu}^2A(z)}}{M_{\nu}\sqrt{A(z)}}=:
	\frac{F_{\nu}(z)}{\sqrt{A(z)}},$$
	where $P_{2\nu}(z):=e^{i\nu\varphi}\tau_{\nu}(\varphi),z=e^{i\varphi},$ has also modulus 1 for 
	$z\in\Gamma_E.$  The function 
	$F_{\nu}$ is algebraic and has no finite poles, it has as branch points 
	$e^{i\varphi_j},j=1,2,3,4,$ only, hence $F_{\nu}(z)=P_1(z)+\sqrt{R(z)}P_2(z),$ 
	where $P_1,P_2$ are polynomials. By \eqref{eq-trig4} and \eqref{trigf} 
	we have $F_{\nu}(z)=P_{2\nu}(z)+\sqrt{R(z)}Q_{2n-2}(z),$ where 
	$Q_{2\nu-2}(z):=e^{i(\nu-1)\varphi}\sigma_{\nu-1}(\varphi).$ Let us normalize the polynomial $Q_{2\nu-2}(z)$ in such a way that 
	$F_{\nu}(z)$ has a pole at $\infty_1,$ where $\infty_1$ is the point infinity 
	in the first sheet of the Riemann surface of the function 
	$w=\sqrt{R(z)}$ associated with $\bar{\mathbb{C}}\setminus \Gamma_E.$ Since 
	the variation of the argument of $F_{\nu}(z)$ when $z$ goes around $\Gamma_{E_j}$ in 
	the clockwise direction is equal to $-2\pi q_j^{(\nu)},$ where $q_j^{(\nu)}$ denotes the number of zeros of $\tau_{\nu}(\varphi)$ on 
	$E_j$,  the total variation of the argument of 
	$F_{\nu}(z)$ when $z$ goes around the boundary of $\bar{\mathbb{C}}\setminus\Gamma_E$ is equal 
	to $-4\pi \nu.$ Hence by the Argument principle $-2\nu=Z-P,$ where $Z,P$ denotes the 
	number of zeros and poles of $F_{\nu}$ in $\bar{\mathbb{C}}\setminus\Gamma_E,$ 
	respectively. Taking into account the choice of the branch of $\sqrt{R}$ we 
	have $P=2\nu,$ hence $Z=0.$ 
	
	Since by \eqref{eq-trig4}
	\begin{equation}
			(P_{2\nu}(z)+\sqrt{R(z)}Q_{2\nu-2}(z))(P_{2\nu}(z)-\sqrt{R(z)}Q_{2\nu-2}(z))
			=M_{\nu}^2A(z)z^{2\nu-a},
		\label{eq-corner2}
	\end{equation}
	we get 
	\begin{equation}
		P_{2\nu}(z_j)=(\sqrt{R}Q_{2\nu-2})(z_j),\quad j=1,\ldots,m^*,
		\label{eq-corner3}
	\end{equation}
	and (b) is proved.

       (b)$\Leftrightarrow$(c). Follows by Lemma 2.	
	
	(d)$\Rightarrow$(c). The proof is analogous to the proof of Corollary 1. One 	applies Theorem 1 with $\mu=1,p=0,W\equiv 1,$ and takes into account the 
	uniqueness of the orthogonal polynomials which have maximal orthogonality (cf. \cite{pehste0}). 
	
	(c)$\Rightarrow$(d) Put $P_{2\nu},Q_{2\nu-2}$ as in the proof of 
	(a)$\Rightarrow$(b). Then we get from 
	\eqref{eq-corner2} and \eqref{eq-corner3} the desired result by applying 
	Theorem 1. The formulas \eqref{eq-trig6},\eqref{eq-trig7} are also obtained by 
	Theorem 1.
	
	(b)$\Rightarrow$(a). The proof is analogous to the proof of \cite[Corollary 
	3.2(a)]{pehste2}. One needs to take into account also the variation of the 
	argument of the function $\mathcal{F}_{\nu}$ from the proof of (a)$\Rightarrow$(b) 
	which can be determined easily by \eqref{eq-trig7}. $\Box$
	
	\newtheorem{remark}{Remark}
	
	\begin{remark} For $\mathcal{A}\equiv 1$ special cases of Corollary 2 can be 
	found in \cite{kru,pet}; the connection with orthogonal polynomials was discovered
       in \cite{pehste2} for any number of arcs even. For $E=[0,2\pi]$ the problem was considered by Szeg\H{o}
	\cite{sze2}.
	\end{remark}
	
\section{Behaviour of zeros}\label{sec.3 }

Let $S_1$ and $S_2$  be neighbourhoods of the arcs $\Gamma_{E_1}$ and  
$\Gamma_{E_2},$ 
respectively. For technical reasons it is more convenient  to take them as 
images of the strips  $S_1^{\Box} $ and $S_2^{\Box} $ under the map 
$x=\phi(u)$, where $$S_1^{\Box}=\left\{-\varepsilon<\Re u<0,-K'<\Im u<K'\right\} $$ 
and
$$S_2^{\Box}=\left\{-K<\Re u<K+\varepsilon,-K'<\Im u<K'\right\}. $$
 
\begin{theorem}\label{thm.2}
	 Let $\mathcal{L} $ be positive definite, $n\geq a+2+v.$ Then   
	  the number of zeros $k_n^{(1)} $ and $k_n^{(2)}$  of the 
	polynomial $P_n$ in $S_1$ and $S_2$ are given for sufficiently large 
	 $n$ by the formulas
	$$k_{n}^{(1)}=n-\frac{1}{2 }(l_n+\beta_n(1-\delta^{(n)}))-
		\sum_{j=1}^{m^* }\frac{1-\lambda_j}{2 }m_j-\frac{1-\delta^{(n)} }{2 }+w-\frac{w_1}{2 }        $$ 
		and$$k_n^{(2)}=\frac{1}{2 }(l_n-\gamma_n(1-\delta^{(n)})-w_2) ,  $$ 
where $\beta_n=1$ if $b^{(n)}\in S_1^{\Box}$ and $\beta_n=0$ if 
$b^{(n)} \notin S_1^{\Box};$  analogously $\gamma_n=1$ if $b^{(n)}\in S_2^{\Box}$ 
and $\gamma_n=0$ if $b^{(n)} \notin S_2^{\Box}.$  Let us point out that $b^{(n)},l_n$ and 
$\delta^{(n)}$ are given uniquely by \addtocounter{z2}{-1}
(\arabic{z2}) and 
\addtocounter{z2}{1}
(\arabic{z2}). 
	\end{theorem}
	
	{\bf Proof.}
	
	First of all let us note that because of the positive definiteness of the 
	functional $\mathcal{L} $  the number $\mu$  from Lemma 2 is equal to $0$. 
	Furthermore let us show that $p=0,$ i.e.,  $P_n(0)\neq 0.$ Indeed, assume that $p=1$ then it follows by \eqref{eq-el10}  that 
	$Q_{n+2-2v}(0)=0.$ Dividing relation \eqref{eq-el10} by $z^2$  one 
	gets that the polynomial $P_{n-1}(z):=P_n(z)/z$ is orthogonal with respect to 
	$\mathcal{L} $ for $j\in(0,n-1]$ and thus $\partial g_{(n)}=0,$ which is by Theorem 1 equivalent to $\delta^{(n)}=0.$ But by Corollary 1 $\delta^{(n)}\in\{-1,1\}$ which is a contradiction. 
	
	Thus the formula \eqref{eq-el15}  can be written as follows:
$$\Psi_n(\phi(u))=ce^{-i\pi m^{(n)}u/K} \left[\frac{H(u+\zeta)}{H(u-\zeta)}\right]^{n+w-a} 
$$
\begin{equation}
\cdot\left[\frac{H(u+\bar{\zeta})}{H(u-\bar{\zeta})}\right]^{n+w-a-1} 
\prod_{j=1}^{m^*}\left[\frac{H(u+v_j)}{H(u-v_j)}\right]^{\lambda_jm_j}\left[\frac{H\left(u+b^{(n)}\right)}
{H(u-b^{(n)})}\right]^{\delta^{(n)}}.	
	\label{eq-42prime}
\end{equation} 
Now we can determine easily all poles of $\Psi_n(\phi(u))$, naturally they 
coincide with the poles of $1+\Psi_n(\phi(u)),$ and are of the same order. In particular
by \eqref{eq-42prime}   $1+\Psi_n(\phi(u))$ 
has $2n+2w$ poles in the parallelogramm of periods $\mathcal{P}. $ 

First let us prove the statement for the case:
\begin{equation}
	\label{case}
	\delta^{(n)}=1\mbox{ and }  W(e^{i\varphi_1})W(e^{i\varphi_2})\neq 0.
	\end{equation}

We suppose that $\varepsilon>0$ is sufficiently small  and such that there are no $v_j$'s 
in $\overline{S_1^{\Box}} $  and in   $\overline{S_2^{\Box}} $  and  
no $b^{(n)}$'s on $\partial S_1^{\Box} $   or  $\partial S_2^{\Box}. $ 

We claim that
\begin{quote}
	{\bf A}: A point $z,|z|<1,$ is a zero of  $P_n$ if and only if it is a zero of the 
	$1+\Psi_n(z).$
	\end{quote}
	
	Let us proof claim A.
From \eqref{eq-el13}  and \eqref{eq-el14}  it follows that 
$$1+\Psi_n(z)=\frac{2P_n(z)(W(z)P_n(z)+\sqrt{R(z)}Q_{n+2-2v}(z))  }{z^{n-(a+1-w)}A(z)g_{(n)}(z)   }. $$ 
 Comparing it with the definition \eqref{eq-el18} of  $\Psi_n(z)$   
it can also be written in the form
\begin{equation}
1+\Psi_n(z)=\frac{2P_n(z)\Psi_n(z)W(z)}{P_n(z)+\sqrt{\frac{V(z)}{W(z) } }Q_{n+2-2v} (z)  }. 	
	\label{proof1}
\end{equation}
Because of the positive definiteness of $\mathcal{L} $ and by 
\cite[Prop.2.3]{pehste0} the polynomials $P_n$ 
and $Q_{n+2-2v}$ have no common zeros, hence by \eqref{proof1}  all
zeros of 
$P_n$ will be zeros of $1+\Psi_n(z).$ 
By the positive definiteness of $\mathcal{L} $ $P_n(z)$ has $n$ zeros in 
$|z|<1.$ Hence since $\phi(u)$ is even  $P_n(\phi(u))$
has $n$ zeros in
$$\Box^+=\{u:-K<\Re u<0,0<\Im u<K'\}$$
and $n$ zeros in $-\Box^+.$ Thus 
$1+\Psi_n(\phi(u))$ has $2n$ zeros at the zeros of  $P_n(\phi(u)).$ 
  Moreover all zeros of  $W(\phi(u)$ will be zeros 
of $1+\Psi_n(\phi(u))$ also. Altogether we found $2n+2w$ zeros of  $1+\Psi_n(\phi(u)).$
By the ellipticity of  $1+\Psi_n(\phi(u))$  the  number 
of zeros and  poles in $\mathcal{P} $ is the same. Since we have shown at the 
beginning of the proof that  $1+\Psi_n(\phi(u))$ has $2n+2w$ poles in $\mathcal{P} $, 
the zeros of  
$P_n(\phi(u))$ and of   $W(\phi(u)$ are the only zeros of  $1+\Psi_n(\phi(u))$ 
in $\mathcal{P}.$ Hence claim {\bf A} is proved. In particular, the number of zeros of  
$P_n(\phi(u))$ and of  $1+\Psi_n(\phi(u))$ in $S_1^{\Box}$ is equal. 
Furthermore, as it is easily seen from \eqref{eq-42prime},  
$1+\Psi_n(\phi(u))$  has one pole in $S_1^{\Box}$, if $b^{(n)}$ is in 
$S_1^{\Box}$ (recall that $\delta^{(n)}=1),$ hence by the Argument
principle    
\begin{equation}
	2\pi(k_n^{1}-\beta_n)=\mbox{var}\arg\limits_{u\in\partial S_1^{\Box} }(
1+\Psi_n(\phi(u))),   
	\label{zeros}
\end{equation} 
where $\partial S_1^{\Box}$ is passed around counterclockwise. Because of the 
ellipticity of the function  $1+\Psi_n(\phi(u))$ we have 
$$\mbox{var}\arg\limits_{u\in\partial S_1^{\Box} }( 1+\Psi_n(\phi(u)))=
	\mbox{var}\arg\limits_{u\in S_1^{(1)} }( 1+\Psi_n(\phi(u)))-$$
\begin{equation}
		-\mbox{var}\arg\limits_{u\in S_1^{(2)} }( 1+\Psi_n(\phi(u)))=:A_1-A_2,
	\label{splitz}
\end{equation} 
where $S_1^{(1)}=\{u:\Re u=0,-K'\leq\Im u\leq K'\}, $  $S_1^{(2)}=\{u:\Re u=-\varepsilon,-K'\leq\Im u\leq K'\}. $ 

To compute $A_1$ we will  describe the range  of the function 
$\Psi_n(\phi(u))$ and  compare it with the range    
of the function $1+\Psi_n(\phi(u)),$ when $u$ varies along $S_1^{(1)}.$ 

For that reason let us write relation \eqref{eq-42prime} in the form  
\begin{equation}
	\Psi_n(\phi(u))=f_n(u)h_n(u),
	\label{split1}
\end{equation} 
where
\begin{equation}
		f_n(u)=c\left(\frac{H(u+\zeta)H(u+\bar{\zeta})}{H(u-\zeta)H(u-\bar{\zeta}) } \right) ^{n-a-1+w} 
		\prod_{j=1}^{m^*}\left[\frac{H(u+\bar{v}_j)}{H(u-v_j)}\right]^{\lambda_jm_j},
	\label{fn1}
\end{equation} 
and
\begin{equation}
	h_n(u)=e^{-i\pi 
	m^{(n)}u/K}\frac{H(u+\zeta)H\left(u+b^{(n)}\right)}{H(u-\zeta)H\left(u-b^{(n)}\right) } 
	\label{hn1}
\end{equation}
$$\cdot\left(-e^{-i\pi u/K} \right)^{\sum_{j\in J}^{ }\lambda_jm_j }\cdot 
e^{-i\pi\sum_{j\in J}^{ }\lambda_j\Re v_j/K }.   $$  
 From Lemma 1 it follows that 
\begin{equation}
	|f_n(u)|=1 \mbox{ for }  u\in S_1^{(1)}, \mbox{ and } |f_n(u)|>1 \mbox{ for }   u\in S_1^{(2)}. 
	\label{fn1prop}
\end{equation} 
Furthermore, for $u\in S_1^{(1)}$  we have by straightforward calculations
\begin{equation}
	|h_n(u)h_n(-u)|=1.
	\label{hn1prop}
\end{equation} 
Recall also (cf.\eqref{eq-el20}) that
\begin{equation}
\Psi_n(\phi(-u))\Psi_n(\phi(u))\equiv 1.	
	\label{psiprop}
\end{equation} 
Now it follows from \eqref{split1}  that
$$\mbox{var}\arg\limits_{u\in S_1^{(1)} }\Psi_n(\phi(u))=
	\mbox{var}\arg\limits_{u\in S_1^{(1)} }f_n(\phi(u))+$$
\begin{equation}
			+\mbox{var}\arg\limits_{u\in S_1^{(1)} }h_n(\phi(u)):=A_{11}+A_{12}.
	\label{split11}
	\end{equation} 
 To compute $A_{11}$ one has to observe firstly, that by Lemma 1 the function
$$\sigma(u,\zeta):=\arg\frac{H(u+\bar{\zeta})}{H(u-\zeta) } $$ 
is equal to the harmonic conjugate of the Green's function 
$g_{\mathbb{C}\setminus\Gamma_E } (\phi(u),\infty)=:g(u)$ hence for $u\in[-iK',iK'] $ 
$$\frac{\partial\sigma}{\partial y } =\frac{\partial g}{\partial x }$$
(by the Cauchy-Riemann conditions, with $u=x+iy),$   and it is obvious that 
$\frac{\partial g}{\partial x }<0$  for $u\in[-iK',iK'] $  and for 
$u\in[-K-iK',-K+iK'] $  we have $\frac{\partial g}{\partial x }>0.$ So 
$\sigma(u,\zeta)$ is strictly decreasing along  $u\in[-iK',iK'] $  and along
$u\in[-K+iK',-K-iK'] ,$  hence
$$\mbox{var}\arg_{ u\in[-iK',iK']	 }\frac{H(u+\bar{\zeta})}{H(u-\zeta) }=
\arg\left(\frac{H(iK'+\bar{\zeta})}{H(iK'-\zeta) }\cdot\frac{H(-iK'-\zeta)}{H(-iK'+\bar{\zeta} ) }\right)+$$
\begin{equation}
		 +2\mu\pi=
-2\pi\Re\zeta/K+2\mu\pi,
	\label{vararg1}
\end{equation} 
where $-\mu\in\mathbb{N}$ and the last equality follows by \eqref{foreta}.  But \eqref{vararg1} holds for any $\zeta\in\Box,$ 
hence by continuity with respect to $\zeta$ for $\Im\zeta=0$    the relation 
\begin{equation}
	-2\pi\zeta/K+2\mu\pi =	\mbox{var}\arg_{ u\in[-iK',iK']	
		 }\frac{H(u+\zeta)}{H(u-\zeta) }
	\label{vararg2}
\end{equation} 
holds. In an analogous way
\begin{equation}
		\mbox{var}\arg_{ u\in[-K+iK',-K-iK']	
		 }\frac{H(u+\bar{\zeta})}{H(u-\zeta) }=2\pi\Re\zeta/K+2\nu\pi,
	\label{vararg3}
\end{equation} 
with $-\nu\in\mathbb{N}_0,$ and for  $\Im\zeta=0$
\begin{equation}
	2\pi\zeta/K+2\nu\pi =	\mbox{var}\arg_{ u\in[-K+iK',-K-iK']	
		 }\frac{H(u+\zeta)}{H(u-\zeta) }.
	\label{vararg4}
\end{equation} 
The variations of the argument of \eqref{vararg2} and  \eqref{vararg4} were 
computed in \cite{pehel2}, but for the sake of completeness let us give another 
proof. Adding \eqref{vararg1} and \eqref{vararg3} one gets
$$	\mbox{var}\arg_{ u\in[-K+iK',-K-iK']\cup[-iK',iK']		
		 }\frac{H(u+\bar{\zeta})}{H(u-\zeta) }=2(\mu+\nu)\pi,$$  
and that the variation of the argument is equal to $-2\pi$ because of the Argument 
principle. Hence $\mu=-1,\nu=0,$ and 
\begin{equation}
		\mbox{var}\arg_{ u\in[-iK',iK']	
		 }\frac{H(u+\bar{\zeta} )}{H(u-\zeta) }=-2\pi\Re\zeta/K-2\pi,
	\label{vararg5}
\end{equation} 
\begin{equation}
	\mbox{var}\arg_{ u\in[-K+iK',-K-iK']	
		 }\frac{H(u+\bar{\zeta} )}{H(u-\zeta) }=2\pi\Re\zeta/K.
	\label{vararg6}
	\end{equation} 
Let us observe,  \eqref{vararg6} with the help of \cite[formula (4.3)]{wid2}  
gives us \eqref{eq-el6}   immediately.
Similarly one computes
$$\mbox{var}\arg_{ u\in[-iK',iK']}\frac{H(u+\bar{v}_j) }{H(u-v_j) }
=-2\pi\Re v_j/K-2\pi,$$
and
$$\mbox{var}\arg_{ u\in[-K+iK',-K-iK']}\frac{H(u+\bar{v}_j) }{H(u-v_j) }
=2\pi\Re v_j/K.$$ Thus we get by \eqref{fn1},\eqref{hn1},\eqref{split11},
\eqref{vararg5} and \eqref{vararg6}     
\begin{equation}
	A_{11}=(n-a-1+w)\left(-2\pi(2\Re\zeta/K+2)\right)-2\pi\sum_{j=1}^{m^* }\lambda_jm_j(\Re v_j/K+1)  
	\label{1arg11}
\end{equation} 
and
\begin{equation}
	A_{12}=-2\pi(\Re\zeta+\Re b^{(n)})/K-4\pi. 
	\label{1arg12}
\end{equation} 
Now we are ready to study the ranges $(\mathcal{C})$ of the function $\Psi_n(\phi(u))$
and to compare it with the range
$(\tilde{\mathcal{C} })  $ of the function $1+\Psi_n(\phi(u))$  when $u$ varies along $S_1^{(1)}.$ 

It was mentioned before that the function
$$\sigma(u,\zeta):=\arg\frac{H(u+\bar{\zeta})}{H(u-\zeta) } $$ 
is strictly decreasing when $u$ varies along $S_1^{(1)}.$ 
The function 
$$\arg\frac{H(u+\zeta)}{H(u-\bar{\zeta}) } $$
 has an analogous property, hence by \eqref{split1},\eqref{fn1} 
for sufficiently large $n$ the argument of  $\Psi_n(\phi(u))$ is strictly 
monotonically decreasing when $u$ varies along $S_1^{(1)}.$ Now using relation \eqref{psiprop} it follows that the curve $\mathcal{C}$ consists 
of two  closed curves with end points 1 (recall the supposition 
$W(e^{i\varphi_1})W(e^{i\varphi_2})\neq 0)$  such that the second one is the 
image of the first one under  reflection with respect to the 
unit circle and to the real axis. Since $\mathcal{C}$ does not run 
through -1, we get  
\begin{equation}
	\mbox{var}\arg_{ u\in S_1^{(1)}=[-iK',iK']}\Psi_n(\phi(u))=2\mbox{var}\arg_{
u\in[-iK',iK']	
		 }(1+\Psi_n(\phi(u)))=2A_1.
	\label{proof2}
\end{equation} 
Next let us show that for $n>N_0$
\begin{equation}
|\Psi_n(\phi(u))|>1 \mbox{ for }u\in 
S_1^{(2)}. 	
	\label{modpsi}
\end{equation} 
Indeed, since by  Lemma 1
$$\left|\frac{H(u+\zeta)H(u+\bar{\zeta})}{H(u-\zeta)H(u-\bar{\zeta}) }\right|>1 
\mbox{ on } u\in 
S_1^{(2)} $$ 
it follows that there exists an $N_0$  such that for any 
$n\in\mathbb{N},n>N_0, $ 
\begin{equation}
	\inf_{u\in S_1^{(2)} } 
	\left|\frac{H(u+\zeta)H(u+\bar{\zeta})}{H(u-\zeta)H(u-\bar{\zeta}) }\right| ^{n-a-1+w} 
	\geq q_{N_0}^n\sup_{u\in S_1^{(2)}}
	\left|\prod_{j=1}^{m^*}\left[\frac{H(u+\bar{v}_j)}{H(u-v_j)}\right]^{\lambda_jm_j}
	\right|\left/
	|h_n(u)|\right.,   
	\label{proofhn}
\end{equation} 
with $q_{N_0}>1, $  wich proves in view of \eqref{split1} and \eqref{fn1} the 
claim.

Thus $$-A_2=\mbox{var}\arg_{u\in S_1^{(2)} }\Psi_n(\phi(u))=-2\pi\beta_n-\mbox{var}\arg_{u\in S_1^{(1)} }\Psi_n(\phi(u))
=-2\pi\beta_n+2\pi l_n,$$
where the first equality follows by the definition \eqref{splitz} of $A_2$ and \eqref{modpsi}, the second 
one uses the ellipticity of   $\Psi_n(\phi(u))$ and the assumption 
$\delta^{(n)}=1$, and the third equality follows by \eqref{split11},\eqref{1arg11},
\eqref{1arg12} and (\arabic{z1}). Hence
\begin{equation}
	A_2=2\pi\beta_n-2\pi l_n.
	\label{1arg2}
\end{equation} 
Finally we get with the help of
\eqref{zeros},\eqref{splitz},\eqref{split11}-\eqref{1arg11},\eqref{proof2}
that 
$$k_n^{(1)}=\frac{1}{K}\left((n-a-1)\Re\zeta+\frac{1}{2 } \sum_{j=1}^{m^* }\lambda_jm_j\Re 
v_j+ \frac{1}{2 }(\Re\zeta+\Re b^{(n)}) \right)+n-a+\frac{1}{2 }\sum_{j=1}^{m^* }\lambda_jm_j  .$$ 
which is the assertion under the assumption (\ref{case}).

 If 
		 $\delta^{(n)}=-1$  and $W(e^{i\varphi_1})W(e^{i\varphi_2})\neq 0$ 
		 then in \eqref{zeros} $\beta_n$ should be omitted, and in  
	\eqref{1arg2} it should appear with minus sign. 	
	
For the calculation of $k_n^{(1)}$ in the case $W(e^{i\varphi_1}) \neq 
0,W(e^{i\varphi_2})=0,$ one  can not repeat the considerations from above without 
any modification since the curve $(\tilde{\mathcal{C} })  $ goes through the 
point 0. Hence one needs to take a 
modified ``interval'' 
$\tilde{J}_{\tilde{\varepsilon} }=[-iK',-i\tilde{\varepsilon}]\cup[i\tilde{\varepsilon},iK']\cup
C^-_{\tilde{\varepsilon} },     $  where $C_{\tilde{\varepsilon} }$  is the 
circumference with the center $u=0$ and the radius $\tilde{\varepsilon},$   and 
$C^-_{\tilde{\varepsilon} },C^+_{\tilde{\varepsilon} }$ 
 are its left- and right-hand halves respectively.
 
 Let $B,D,O$ denote the images   of the 
 points $-i\tilde{\varepsilon},i\tilde{\varepsilon},0$ under the function 
 $\Psi_n(\phi(u)).$  Note $O$ is the 
 point $-1.$ Since the variation of the argument of $\Psi_n(\phi(u))$ 
 (for $n$ large enough) is strictly decreasing along $[-i\tilde{\varepsilon},
 i\tilde{\varepsilon}], $  the curve $BOD$ is such that $\frac{\pi}{2 }<\arg 
 D<\pi,\pi<\arg B<\frac{3\pi}{2 }.  $  Now the variation of the argument of the 
 function $1+\Psi_n(\phi(u))$  along the circumference $C_{\tilde{\varepsilon} }$  
 is equal to $2\pi$ (recall that $\Psi_n(\phi(0))=-1$ since $W(\phi(0))\neq 0).$ Thus the image of $C_{\tilde{\varepsilon} }$  under the 
 function $\Psi_n(\phi(u))$  is such that the point $-1$ lies inside of
 $\Psi_n(\phi(C_{\tilde{\varepsilon} }))$  and the curve $\Psi_n(\phi(u))$ goes 
 counterclockwise around $O$ when $u$ varies counterclockwise around $u=0$ 
 along $C_{\tilde{\varepsilon} }.$
 Hence the image of $C^-_{\tilde{\varepsilon} }$  will  be such that $O$ is at 
 the 
 right hand side of $\Psi_n(\phi(C^-_{\tilde{\varepsilon} })).$ Finally we get that  
 the variation of the argument of $\Psi_n(\phi(u))$ along $[0,iK'] $  is equal 
 to $-2(2\kappa+1)\pi,\kappa\in\mathbb{N}_0, $ because of $\Psi_n(\phi(0))=-1$ 
 and $\Psi_n(\phi(iK'))=1.$  Thus the   
 considerations give finally
 $$\mbox{var}\arg_{ u\in\tilde{J}_{\tilde{\varepsilon} }	
		 }(1+\Psi_n(\phi(u)))=-2(\kappa+1)\pi.$$ 
	
		 The variation of the argument of the function $1+\Psi_n(\phi(u))$  along 
$\partial S_2^{\Box}$  is calculated in an analogous way.
		 Other cases are considered in the same manner.      $\Box$ 
		 
	\begin{theorem}\label{thm.3}
	Let the functional $\mathcal{L} $  be positive definite,  let $z_{j,n},j=1,\ldots,n, $ 
	be the zeros of $P_n$  and let $Z$ be the set of all accumulation points of 
	$(z_{j,n})^{n,\infty}_{j=1,n=1}.   $ Furthermore put
\begin{equation}
S=\phi(L), L=\{u\in\Box:\Im u=\Im\zeta+K'\}, 
\Xi=\{e^{i\xi_j} :\lambda_j=-1\}.	
	\label{curve}
\end{equation}
 Then the following statements hold:
	\begin{alphenumerate}
	\item
 $Z\subseteq\Gamma_E\cup \Xi \cup S,$ 
where $Z\cap\Gamma_E=\Gamma_E$ and $Z\cap\Xi=\Xi.$
\item 
 $Z\cap S=S$ and thus $Z=\Gamma_E\cup \Xi \cup S$ if the harmonic measure $\omega_2(\infty)$ 
of  $\Gamma_{E_2}$ is an irrational number.
\item
If  
$\omega_2(\infty)$  is  rational, then $\mathcal{N}:=Z\cap S$  is a
finite set and $\phi^{-1}(\mathcal{N})=$
$$=\left\{\left(\Re\zeta(2n+2w-2a-1)+\sum_{j=1}^{m^* }\lambda_jm_j
\Re v_j \right)/K+i(\Im\zeta+K'):n\in\mathbb{N}  \right\} \cap\Box  $$ 
\end{alphenumerate}
\end{theorem}
{\bf Proof.}

First let us recall (see claim A in the proof of Theorem 2) that $P_n(\phi(u))$ has a zero at 
$u\in\inte\Box$ if and only if $\Psi_n(\phi(u))=-1.$  With the 
help of relation \eqref{modpsi} and by the continuity of $\Psi_n(\phi(u))$ it 
follows by representation
\eqref{eq-42prime} that each mass-point $\phi(v_j)$ of $\mathcal{L}$ is an accumulation 
point of zeros of $\left(P_n(\phi(u))\right)$ since 
$\lambda_j=-1$ and thus $v_j$ is a zero of $\Psi_n(\phi(u)).$  Furthermore, it follows by the 
same reasons that other accumulation points of zeros of $(P_n)$ in 
$\mathbb{C}\setminus E,$ more precisely in $\{|z|\leq 1\}\setminus E,$ since $P_n(z)$ 
has all zeros in $\{|z|<1\},$ may appear at accumulation points of 
$(\phi(b^{(n)})),b^{(n)}\in \inte\Box_+,$ only, where an accumulation point 
$b^*\in\bar{\Box}_+$ of $(b^{(n)})$, i.e., $b^*=\lim_{k\to\infty}b^{(n_k)},$ is 
a limit point of zeros of $\left(P_{n_k}(\phi(u))\right)$ if and only if
$\delta^{(n_k)}=-1$ for $k\geq k_0.$ Next let us show that $b^*\in L.$
Indeed, putting
 $\iota=m^{(n)}+\sum_{j=1}^{m^* }\lambda_jm_j\Im v_j/K',$
  it follows from  
	\addtocounter{z2}{-1}
	(\arabic{z2}) that there are two possibilities for $\iota:$ 
	either $\iota=1$ or $\iota=0.$ Furthermore, for $\delta^{(n)}=-1 $ 
\begin{equation}
	\label{th4.1}
	\Im b^{(n)}=\Im\zeta+K' \mbox{ for  }\iota=1 ,
	\mbox{ and } \Im b^{(n)}=\Im\zeta 
	  \mbox{ for  }\iota=0, 
	\end{equation}  
	and for $\delta^{(n)}=1 $ 
	\begin{equation}
	\Im b^{(n)}=-\Im\zeta-K' \mbox{ for  }\iota=1,
	\mbox{ and } \Im b^{(n)}=-\Im\zeta  
	 \mbox{ for  }\iota=0.	
		\label{th4.2}
	\end{equation} 
	Now from $ \lim_{k\to\infty}b^{(n_k)}= b^*\in\Box_+$ and  $\delta^{(n_k)}=-1$ for $k\geq k_0$ 
	we obtain that the first relation in \eqref{th4.1}, i.e.,  $\iota=1$ holds.Thus what remains to
be shown is the relation $Z \cap \Gamma _E = \Gamma _E$ which follows
by (3). 
       
 Let us mention that it is known even, see e.g.	 [9,Thm.9.2;26], that each isolated mass point attracts exactly one zero of $P_n.$ This can be proved also in the following way: taking into consideration the facts that $m_j=1,$ since 
$\mathcal{L}$ is positive definite, and that there are no $b^{(n)}$'s in the 
neighbourhood of mass points $v_j$ for sufficiently large $n,$ we have
\begin{equation}
	\mbox{var}\arg_{ u\in\mathcal{B} }\Psi_n(\phi(u))=\mbox{var}\arg_{
u\in\mathcal{B}	
		 }(1+\Psi_n(\phi(u))).
	\label{th4.3}
\end{equation} 
But the left hand side expression in \eqref{th4.3} is equal to $2\pi i$ because of 
\eqref{eq-el15}, hence the number of zeros of $P_n$ in the neighbourhood is 
equal to 1 proving the assertion.

	Concerning part b), let us recall Chebyshev's theorem:
	If $\alpha,0<\alpha<1,$ is an irrational number, then for any $x\in\mathbb{R}$ and 
	for any $\varepsilon>0$ it is possible to find $n\in\mathbb{N}$ and  
	$m\in\mathbb{Z}$ such that
	\begin{equation}
		|n\alpha-m-x|<\varepsilon.
		\label{cheb}
	\end{equation} 
	Put $\alpha=-\frac{\Re\zeta}{K}$ and
	$$x=(w-a-1/2)\frac{\Re\zeta}{K}-\sum_{j=1}^{m^* }\lambda_jm_j\frac{\Re v_j}{2K}-
	\frac{\delta b}{2K  } -w_2/2,$$
	where $\delta\in\{-1,1\}$ and $b,-K<b<0,$ are arbitrary. Then for any 
	$\varepsilon>0$ it is possible to find  $n\in\mathbb{N}$ and  
	$m\in\mathbb{Z}$ such that \eqref{cheb} holds. By (\arabic{z2})
	$$n\alpha=(w-a-1/2)\frac{\Re\zeta}{K}+\sum_{j=1}^{m^* }
	\lambda_jm_j\frac{\Re v_j}{2K}+\frac{\delta^{(n)} \Re 
	b^{(n)} }{2K  }+w_2/2-(l_n-w_2)/2,$$
	Inserting this in  \eqref{cheb}  gives
	$$\left|\frac{\delta b}{2K  }-\frac{\delta^{(n)} \Re 
	b^{(n)} }{2K  }-\left((l_n-w_2)/2-m\right) \right| <\varepsilon.$$
	Since both $\frac{\delta b}{2K  }$ and $\frac{\delta^{(n)} \Re 
	b^{(n)} }{2K  }$ are in $(-1/2,1/2)$, we have $(l_n-w_2)/2=m,$ and
	$$\left|\frac{\delta b}{2K  }-\frac{\delta^{(n)} \Re 
	b^{(n)} }{2K  } \right| <\varepsilon.$$
Hence for any $b\in(-K,0)$ and $\delta=-1$ there is a subsequence $(n_k)$ of 
the natural numbers such that $b^{(n_k)}$ satisfying (\arabic{z2}) and \addtocounter{z2}{-1}
	(\arabic{z2}) with $\delta^{(n_k)}=-1$ (recall \eqref{th4.1} and the fact that 
	$\iota=1)$ tends to $b+i(\Im\zeta+K').$ Hence b) is proved.
	
Part c) is proved  in the same manner by taking into account that 
	 there exists only a finite number of possible solutions of 
	\addtocounter{z2}{1}(\arabic{z2}) for all $n\in\mathbb{N}$ and that only for $\delta^{(n)}=-1$ these 
	solutions will attract zeros of $P_n$.
	 $\Box$

	\begin{remark}	 
	  Part c) can be  proved (with description of the set $\mathcal{N}$ in other 
	  terms) by combination of \cite[Thm.3.3]{pehste3},\cite[Remark 
3.1,Thm.4.2]{pehste2} and the
	  corrected version of \cite[Thm.4.2]{pehste4} (compare also \cite[Thm.2]{barlop},\cite{luk1}).
	\end{remark}						 

\begin{remark}
Let us note that Tomchuk \cite{tom} on page 2 before Thm.2 claims that the function at the right hand side in \eqref{eq-cor5}, denoted by him by $p(z,\sqrt{R(z)}),$ has all zeros in $|z|<1.$ We would like to mention that the claim is not correct. Indeed let us assume that the claim is correct. Then the function $\Omega_n(u)$ from \eqref{eq-cor5} and \eqref{eq-cor2} (see the introduction of $\phi$ at the end of Section 2) has all zeros in the upper half of $\Box$, i.e. in $\Box_+=(-K,0)\times(0,iK')$ or in $-\Box_+.$ Moreover by \eqref{eq-cor2} $-\delta^{(n)}b^{(n)}\in\pm\Box_+$ for $n\geq a+1+v.$ But let us show that this is impossible, because there always exist a subsequence $(n_k)$ such that $\delta^{(n_k)}=1$ and $b^{(n_k)}\in[-K,0]\times[-iK',0]$ if the functional $\mathcal L$ is positive definite. Indeed, by 
\addtocounter{z2}{1}(\arabic{z2}) we can choose a sequence $(n_k)$ 
such that for any 
$k\in\mathbb{N}$ $\delta^{(n_k)}=1.$ Now it follows from \eqref{eq-el19}  and 
from \cite[Ths. 2.1,2.2]{pehste1} that
\begin{equation}
	\Psi_{1,n}(z)=\frac{P_n^2(z)\left(F(z)+\frac{\Omega_n(z)}{P_n(z) } \right)^2A(z)V(z) }
	{z^{n+a+1-w}g_{(n)}(z)   },  
	\label{th4.4}
\end{equation}  
where $$F(z)=\mathcal{L}\left(\frac{x+z}{x-z  };{\cal A,W},\lambda) \right) $$
is the Caratheodory function associated with the functional ${\cal L}$ and 
$\Omega_n(z)$ (do not mix with $\Omega_n$ from Corollary 1) are the polynomials 
of second kind. Since $\Psi_{1,n}(\phi(u))=\Psi_n(\phi(-u))$ the functions 
$\Psi_{1,n_k}(\phi(u))$ have a zero at $u=b^{(n_k)}$ by \eqref{eq-el15}. But 
$g_{(n_k)}(\phi(u))$ also has a zero at $u=b^{(n_k)}$, hence by \eqref{th4.4}
$$F(\phi(u))+\frac{\Omega_{n_k}(\phi(u)) }{P_{n_k}(\phi(u))  }$$
has a zero   at $u=b^{(n_k)}.$ Now by \cite[Ths. 12.1,12.2]{ger} the function  
$F(z)+\frac{\Omega_n(z)}{P_n(z) }$ has no zeros inside the unit circle, hence 
$|\phi(b^{(n_k)})|\geq 1 $ and $b^{(n_k)}\in[-K,0]\times[-iK',0].$
	\end{remark}

 
%
%

%
%
\end{document}